\documentclass[11pt]{article}
\usepackage[margin=1in]{geometry}

\usepackage{amsmath,latexsym,amsfonts,amssymb}
\usepackage{mathrsfs,enumerate}
\usepackage{graphicx,epsfig,color}
\usepackage{tikz}
\usepackage{xcolor}
\usepackage{pgfplots}
\usetikzlibrary{arrows.meta,calc,decorations.pathreplacing}
\newcommand{\un}{\underline}

\newcommand{\q}{\mathbf q}
\newcommand{\p}{\mathbf p}

\newcommand{\e}{\mathbf e}

\newcommand{\cI}{\mathcal I}

\newcommand{\cM}{\mathcal M}

\newcommand{\cO}{\mathcal O}
\newcommand{\cP}{\mathcal P}

\newcommand{\bbM}{\mathbb M}
\newcommand{\bbN}{\mathbb N}

\newcommand{\bbP}{\mathbb P}

\newcommand{\bbR}{\mathbb R}

\newcommand{\bbZ}{\mathbb Z}

\newcommand{\Slab}[3]{%
  \coordinate (O)   at #1;
  \coordinate (B)   at ($ (O) + (#2,0) $);
  \coordinate (Od)  at ($ (O) + (0.9,0.45) $);
  \coordinate (Bd)  at ($ (B) + (0.9,0.45) $);
  \coordinate (Oz)  at ($ (O) + (0,0.9) $);
  \coordinate (Bz)  at ($ (B) + (0,0.9) $);
  \coordinate (Odz) at ($ (Od) + (0,0.9) $);
  \coordinate (Bdz) at ($ (Bd) + (0,0.9) $);

  \fill[#3!35] (Odz) -- (Bdz) -- (Bz) -- (Oz) -- cycle; 
  \fill[#3!55] (B) -- (Bd) -- (Bdz) -- (Bz) -- cycle;    
  \fill[#3!25] (O) -- (B) -- (Bz) -- (Oz) -- cycle;      

  \draw[black,line width=.3pt]
    (O) -- (B) -- (Bd) -- (Od) -- cycle
    (Oz) -- (Bz) -- (Bdz) -- (Odz) -- cycle
    (O) -- (Oz) (B) -- (Bz) (Bd) -- (Bdz) (Od) -- (Odz);
}

\title{\bf Multiscale Analysis of the Conductivity in the Lorentz Mirrors Model}

\author{Rapha\"el Lefevere\\[1mm]
\small Laboratoire de Probabilit\'es, Statistiques et Mod\'elisation,\\
\small Universit\'e Paris Cit\'e}

\date{\today}

\begin{document}
\maketitle

\begin{abstract}
We study deterministic transport in a random medium using the mirrors model, 
a lattice Lorentz gas at unit density in which a particle moves deterministically 
in a frozen random configuration of mirror--scatterers. Despite the absence of chaos 
and the existence of infinitely many finite trapping loops, numerical evidence suggests 
that the model exhibits normal conductivity in $d=3$. We develop a recursive multiscale expansion for the crossing probability of a slab of width $N$, showing that 
$C_N \sim \kappa/(\kappa+N)$ for large $N$ and computing the conductivity constant 
$\kappa$ through a renormalization procedure based on scale concatenation.

The key idea is that a slab of width $2^{n+1}$ may be decomposed into two independent
slabs of width $2^n$, and that the crossing event can be expressed as a sum over
trajectories that cross the left half, return a few times to the interface between the
two halves, and finally cross the right half. This gives rise to a recursive relation in which the leading normal conduction scaling is
propagated from scale $2^n$ to scale $2^{n+1}$, while the subleading term yields
the correction to the conductivity constant. For the mirrors model, it involves a correction factor with respect to a reference Markovian process that 
encodes correlations between crossings at scale $2^n$. These correlations are 
controlled through a closure assumption whose structure is shaped by the hard--core exclusion inherent to the reversible deterministic dynamics. The dominant contribution comes from trajectories with a single interface return, whose
asymptotic contribution we identify explicitly.

In $d=3$ the recursion takes the form 
\[
\kappa_{n+1}=\kappa_n\Bigl(1+\alpha\,\frac{\kappa_n}{2^n}+o(2^{-n})\Bigr),
\qquad \alpha\simeq 0.0374,
\]
leading to a finite limit $\kappa_\infty\simeq 1.5403$, in good agreement with 
numerical simulations performed in this paper. This value is close to the conductivity constant of a non--backtracking 
random walk, suggesting that the large--scale behavior of the mirrors model is 
effectively Markovian even though the microscopic dynamics is fully deterministic.

\end{abstract}


\section{Introduction}

Transport in deterministic systems with quenched disorder plays a central role in 
statistical mechanics. A fundamental question is whether a macroscopic law such as 
Fick's law can emerge from a microscopic dynamics that is neither stochastic nor chaotic. 
A paradigmatic example is the lattice Lorentz gas, where a particle moves deterministically 
in a fixed random configuration of scatterers. Even in this simple setting, understanding 
the emergence of diffusion and normal conductivity remains a major challenge.

The mirrors model, introduced in \cite{Ruijgrok}, provides a particularly striking 
instance of such a system. The scatterers are ``mirrors'' placed with a density parameter $0<p\leq 1$ at the sites of 
$\mathbb{Z}^d$, each implementing a random reflection chosen uniformly among those satisfying 
reversibility and the no--U--turn constraint. The resulting dynamics is fully deterministic, 
time--reversible, and non--chaotic; in particular, infinitely many finite loops occur 
with positive probability. This leads to strong memory effects and non--Markovian behavior 
that severely complicate probabilistic analysis. It is straightforward to see that $\bigl(\tfrac{p}{2d - 1}\bigr)^4$ provides a lower bound on the probability that a particle starting from the origin becomes trapped in a finite loop forever. Proving true diffusivity of the motion in this model is therefore an elusive challenge. The best results in this direction are due to Elboim, Gloria, and Hern\'andez  \cite{ElboimGloriaHernandez}, who showed that for all dimensions $d \ge 4$ and sufficiently small $p > 0$, particle trajectories exhibit diffusive behavior up to times of order $1/p^{\alpha}$, with $\alpha > 1$; in particular, trajectories remain open up to this time scale. In contrast, in dimension $d = 2$ with $p = 1$, it is known \cite{Grimmett,bunitroub} that every trajectory eventually closes with probability one, while the case $p<1$ was analyzed in \cite{Sanders}.
In the present work we do not address diffusion or recurrence on $\mathbb{Z}^d$ itself but instead focus on establishing \emph{normal conductivity} in the mirrors model. This property has been established in a dilute regime for a non-lattice random Lorentz gas~\cite{Basile}, and in an anisotropic random lattice Lorentz gas at full density for $d \ge 7$~\cite{Lefevere}. For this purpose, the most fundamental observable is the \emph{crossing probability} of a finite slab~\cite{Lefevere,ChiffaudelLefevere}. Despite the seemingly unfavorable features mentioned above from the viewpoint of diffusion, numerical simulations~\cite{ChiffaudelLefevere} for the case of full density ($p=1$) in $d=3$ support the conjecture that the mirrors model nevertheless displays \emph{normal conductive behavior}. Specifically, the crossing probability $C_N$ of a slab of width $N$ appears in~\cite{ChiffaudelLefevere} to scale as $C_N \sim \kappa / N$, with $\kappa \simeq 1.535$ \footnote{We will see below that the value of the conductivity is updated to $1.5403$, a value still within the confidence interval reported in \cite{ChiffaudelLefevere}.}, a value close to $3/2$. This latter value is the conductivity constant of a non-backtracking random walk, suggesting that despite the strong memory effects inherent to the mirrors dynamics, its large-scale behavior remains close to that of a Markovian process.
This raises two natural questions: (i) how can a deterministic, non--chaotic system 
effectively behave as a Markov process at large scales? and (ii) can the conductivity 
be computed analytically?

The goal of this paper is to answer both questions through a multiscale recursive analysis of the 
crossing probability of slabs. For a slab $\Lambda_N$ of width $N$, we denote by $C_N$ 
the probability (with respect to the disorder) that a particle entering on the left 
exits on the right. The conductivity at scale $N$ is defined as
\[
\kappa_N = \frac{N C_N}{1-C_N},
\]
and our aim is to understand the asymptotic behavior of $\kappa_N$ as $N\to\infty$.
Our strategy is recursive. At scale $2^n$, we assume that the slab crossing
probability already has the normal conduction order of magnitude
\[
c_n=C_{2^n}\asymp 2^{-n},
\]
or equivalently that the effective conductivity
\[
\kappa_n=\frac{2^n c_n}{1-c_n}
\]
remains of order one. The purpose of the multiscale decomposition is then not to
rederive this scaling from scratch, but to compute how the conductivity constant
is corrected when one passes from scale $2^n$ to scale $2^{n+1}$. In other words,
the recursive step takes as input a normal conductive behavior at scale $2^n$ and outputs
a renormalized value of $\kappa_{n+1}$. The main point is that, once the leading
$2^{-n}$ scaling is factored out, the remaining information is contained in the
difference between $\kappa_{n+1}$ and $\kappa_n$, which can be expressed in terms
of transfer correlations inside the two half-slabs.

Our analysis is based on a decomposition of a slab of width $2^{n+1}$ into two slabs 
of width $2^n$. For a non--backtracking random walk, independence yields an exact 
recursion for the crossing probability, from which one derives 
$\kappa_N \equiv d/(d-1)$. In the mirrors model, the left and right halves are 
independent, but a trajectory that repeatedly crosses the interface forces different 
crossing segments to satisfy deterministic compatibility constraints inside each half. 
These induced correlations between successive crossings form the main difficulty. 
A closure assumption is used to control these correlations. Under this decomposition, the 
correlation functions split naturally into two sectors: one in which correlations vanish 
identically due to hard--core constraints imposed by the reversible dynamics, and another 
in which one assumes increasing independence between trajectories as the scale grows.

We introduce quantities $\eta_n(l)$ that measure the deviation from independence for $l$ transfer pieces of a trajectory inside a half-slab of width $2^n$. The conductivity 
correction at scale $2^{n+1}$ is governed by a weighted sum of the $\eta_n(l)$, with 
weights decaying geometrically in $l$. The dominant contribution arises from $l=2$. 
We show that in $d=3$ this contribution converges rapidly to a constant that we compute 
numerically and bound analytically. Thus the role of the quantities $\eta_n(l)$ is to measure, under the inductive
conductive scaling $c_n\asymp 2^{-n}$, the deviation from the Markovian recursion
and hence the correction to the conductivity constant from one scale to the next. The resulting recursion is
\[
\kappa_{n+1}
= \kappa_n\Bigl( 1 + \frac{\kappa_n}{2^n}(\alpha+o(1))\Bigr),
\qquad \alpha\simeq 0.0374.
\]
This leads to a finite limit $\kappa_\infty\simeq 1.5403$, in good agreement with the
simulations presented below and very close to the value $3/2$ for the non--backtracking
random walk.

The restriction to the three-dimensional case is motivated by the relation between the
present slab geometry and the more standard cubic-box geometry considered in
\cite{ChiffaudelLefevere}. In that work, we showed that a necessary condition for normal
conductivity is that $N C_N$ converges to a finite nonzero constant when $N\to\infty$ and
$C_N$ is measured in a cubic box of side length $N$. Numerical results indicate that this
condition behaves very differently in dimensions $2$ and $3$. In dimension $3$, for cubic
boxes, $N C_N$ appears to converge to a constant, and the limiting conductivity agrees
numerically with the one obtained in the slab geometry considered here. By contrast, in
dimension $2$, unpublished numerical work \cite{Chiffaudel} shows that in square boxes the quantity $N C_N$
does not converge to a single constant, but instead oscillates between two values. This is
the reason why we focus here on the case $d=3$, where the slab geometry appears to capture
the same large-scale conductivity as the cubic geometry. One may nevertheless study the
two-dimensional model in the present slab geometry with large transverse section: preliminary
numerics suggest that a limit $\kappa_\infty$ may also exist there. However, in contrast
with the three-dimensional case, that slab limit would not reflect the behavior in square
boxes. For this reason, the two-dimensional case is of a different nature and is not pursued
in the present paper.

\vspace{1ex}
\noindent\textbf{Outline of the paper.}
Section~2 defines the model and the crossing probability. 
Section~3 recalls the Markovian benchmark given by the non--backtracking random walk.  
Section~4 presents the multiscale recursion for the mirrors model.  
Section~5 contains the closure assumption and the derivation of the conductivity recursion.  
Section~6 analyzes the dominant second--order term and computes the limit conductivity.  
Appendix~A contains the computation of $C_1$, and Appendices~B and~C summarize the 
correlation structure and bounds on the error terms.

\begin{figure}[htbp]
\centering
\begin{tikzpicture}[scale=0.8, >=Latex]

  \def\Nx{6}
  \def\Ny{5}

  \foreach \x in {0,...,\Nx} {
    \draw[gray!40] (\x,0) -- (\x,\Ny);
  }
  \foreach \y in {0,...,\Ny} {
    \draw[gray!40] (0,\y) -- (\Nx,\y);
  }

  \foreach \x in {0,...,\Nx} {
    \foreach \y in {0,...,\Ny} {
      \fill[gray!70] (\x,\y) circle (1pt);
    }
  }

  \def\mirrorNE#1#2{%
    \draw[red,very thick] (#1-0.25,#2-0.25) -- (#1+0.25,#2+0.25);%
  }
  \def\mirrorNW#1#2{%
    \draw[red,very thick] (#1-0.25,#2+0.25) -- (#1+0.25,#2-0.25);%
  }

  \mirrorNE{2}{2}
  \mirrorNE{2}{0}
  \mirrorNE{4}{0}
  \mirrorNE{4}{3}

  \mirrorNW{1}{1}
  \mirrorNE{3}{1}
  \mirrorNW{5}{1}
  \mirrorNW{1}{4}
  \mirrorNE{3}{4}
  \mirrorNW{5}{4}
  \mirrorNW{0}{1}
  \mirrorNE{0}{4}
  \mirrorNW{6}{2}
  \mirrorNE{6}{4}
  \mirrorNW{1}{5}
  \mirrorNE{5}{5}


  \draw[->,thick,blue] (-0.7,2) -- (0,2);

  \foreach \p in {(0,2),(1,2),(2,2),(2,3),(2,4),(2,5),
                  (2,0),(3,0),(4,0),(4,1),(4,2),(4,3),(5,3),(6,3)} {
    \fill[blue] \p circle (2pt);
  }

  \draw[->,thick,blue] (0,2) -- (1,2);
  \draw[->,thick,blue] (1,2) -- (2,2);

  \draw[->,thick,blue] (2,2) -- (2,3);
  \draw[->,thick,blue] (2,3) -- (2,4);
  \draw[->,thick,blue] (2,4) -- (2,5);


  \draw[->,thick,blue] (2,5) -- (2,5.6);

  \draw[->,thick,blue] (2,-0.6) -- (2,0);


  \draw[->,thick,blue] (2,0) -- (3,0);

  \draw[->,thick,blue] (3,0) -- (4,0);
  \draw[->,thick,blue] (4,0) -- (4,1);

  \draw[->,thick,blue] (4,1) -- (4,2);
  \draw[->,thick,blue] (4,2) -- (4,3);

  \draw[->,thick,blue] (4,3) -- (5,3);
  \draw[->,thick,blue] (5,3) -- (6,3);

  \draw[->,thick,blue] (6,3) -- (6.7,3);

\end{tikzpicture}

\caption{Two-dimensional mirrors model on a portion of the lattice $\mathbb{Z}^2$.
Mirrors (red) reflect the deterministic trajectory (blue). Periodic vertical
boundary conditions are indicated by short arrows leaving and re-entering
the domain.}
\label{fig:mirrors_2D_clean_correct}
\end{figure}

\section{Model and Definitions}

We consider a subset of $\bbZ^d$ of the following form:
\[
\Lambda_N = \{1,\ldots,M\}^{d-1}\times\{1,\ldots,N\},
\]
with $M$ much larger than $N$, and with periodic or sufficiently large transverse 
extent so that boundary effects in directions orthogonal to $\mathbf e_1$ may be neglected.

Let $(\e_1,\ldots,\e_d)$ be the canonical basis of $\bbR^d$, and define the set of 
possible velocities
\[
\cP = \{\pm\e_1,\ldots,\pm \e_d\}.
\]
The motion of particles, travelling on the edges of $\bbZ^d$ with unit velocity, is 
described at every time $t\in\bbN$ by a point in the phase space
\[
\cM_N = \{(\q,\p): \q\in \Lambda_N,\ \p\in\cP\}.
\]
\begin{figure}[htbp]
\centering

\begin{minipage}{0.47\textwidth}
\centering
\begin{tikzpicture}[scale=0.8, >=Latex]

  \def\Nx{5} 
  \def\Ny{4} 

  \foreach \x in {0,...,6} {
    \draw[gray!40] (\x,0) -- (\x,5);
  }
  \foreach \y in {0,...,5} {
    \draw[gray!40] (0,\y) -- (6,\y);
  }

  \draw[thick] (1,1) rectangle (\Nx,\Ny);
  \node at (\Nx/2,\Ny/2) {$\Lambda_N$};

  \draw[thick] (1,1) -- (1,\Ny);
  \draw[thick] (\Nx,1) -- (\Nx,\Ny);

  \node[below] at (1,0) {$q_1=1$};
  \node[below] at (\Nx,0) {$q_1=N$};

  \foreach \y in {1,...,\Ny} {
    \fill (1,\y) circle (2pt);                    
    \draw[->,thick,blue] (1,\y) -- (1+0.7,\y);    
  }
  \node[above] at (1,\Ny+0.4) {$I_N^+$};

  \foreach \y in {1,...,\Ny} {
    \fill (\Nx,\y) circle (2pt);                  
    \draw[->,thick,blue] (\Nx,\y) -- (\Nx-0.7,\y);
  }
  \node[above] at (\Nx,\Ny+0.4) {$I_N^-$};

\end{tikzpicture}

\vspace{1ex}
{\small Incoming phase-space boundaries $I_N^+$ and $I_N^-$.}

\end{minipage}
\hfill
\begin{minipage}{0.47\textwidth}
\centering
\begin{tikzpicture}[scale=0.8, >=Latex]

  \def\Nx{5}
  \def\Ny{4}

  \foreach \x in {0,...,6} {
    \draw[gray!40] (\x,0) -- (\x,5);
  }
  \foreach \y in {0,...,5} {
    \draw[gray!40] (0,\y) -- (6,\y);
  }

  \draw[thick] (1,1) rectangle (\Nx,\Ny);
  \node at (\Nx/2,\Ny/2) {$\Lambda_N$};

  \draw[thick] (1,1) -- (1,\Ny);
  \draw[thick] (\Nx,1) -- (\Nx,\Ny);

  \node[below] at (1,0) {$q_1=1$};
  \node[below] at (\Nx,0) {$q_1=N$};

  \foreach \y in {1,...,\Ny} {
    \fill (0,\y) circle (2pt);                    
    \draw[->,thick,blue] (0,\y) -- (-0.7,\y);
  }
  \node[above] at (0,\Ny+0.4) {$O_N^-$};

  \foreach \y in {1,...,\Ny} {
    \fill (\Nx+1,\y) circle (2pt);                
    \draw[->,thick,blue] (\Nx+1,\y) -- (\Nx+1+0.7,\y);
  }
  \node[above] at (\Nx+1,\Ny+0.4) {$O_N^+$};

\end{tikzpicture}

\vspace{1ex}
{\small Outgoing phase-space boundaries $O_N^-$ and $O_N^+$.}

\end{minipage}

\caption{The slab $\Lambda_N$ as a subset of $\mathbb Z^d$ (here represented in two 
dimensions), together with the incoming sets $I_N^\pm$ (left) and outgoing sets 
$O_N^\pm$ (right) in the direction $\mathbf e_1$. The underlying lattice structure 
is shown explicitly as the edges of $\mathbb Z^2$, and the thick dots mark the 
phase-space positions where the velocities are attached.}
\label{fig:LambdaN_IN_ON_lattice}
\end{figure}
We define the outgoing boundaries corresponding to particles leaving the slab on the 
left or right with appropriate velocities:
\begin{align}
O^-_N &= \{(\q-\e_1,-\e_1): q_1=1\},\\
O^+_N &= \{(\q+\e_1,\e_1): q_1=N\},
\end{align}
and $O_N = O_N^+\cup O_N^-$. Similarly, the incoming sets are
\begin{align}
I^-_N &= \{(\q,-\e_1): q_1=N\},\\
I^+_N &= \{(\q,\e_1): q_1=1\},
\end{align}
with $I_N = I_N^-\cup I_N^+$.

On the vertices of $\Lambda_N$, we put ``mirrors'' that modify the orientation of 
the velocity of the particles. For each $\q\in \Lambda_N$, the action of a mirror on 
the velocity of an incoming particle is represented by a bijection $\pi(\q;\cdot)$ 
of $\cP$ into itself. It satisfies the reversibility condition:
\[
\pi(\q;-\pi(\q;\p))=-\p,\qquad \forall\, (\q,\p)\in\cM_N,
\]
and the no U--turn property:
\[
\pi(\q;\p)\neq -\p.
\]

The dynamics is defined on $\cM_N$ by
\begin{equation}
F(\q,\p;\pi) = \bigl(\q+\pi(\q;\p),\,\pi(\q;\p)\bigr).
\label{eq:Fdef}
\end{equation}
It is straightforward to check that the map $F$ is bijective between $\cM_N$ and 
$(\cM_N\backslash I_N)\cup O_N$.
Indeed, $F$ is injective because if
\[
F(\q,\p;\pi)=F(\widetilde \q,\widetilde \p;\pi)=(\q',\p'),
\]
then $\q=\q'-\p'=\widetilde \q$, and since $\pi(\q;\cdot)$ is a bijection we also get
$\p=\widetilde \p$. It is surjective onto $(\cM_N\setminus I_N)\cup O_N$ because for every
$(\q',\p')$ in that set, the candidate preimage must be $(\q'-\p',\p)$ with
$\pi(\q'-\p';\p)=\p'$, and such a unique $\p$ exists since $\pi(\q'-\p';\cdot)$ is a
bijection.

At each site of $\bbZ^d$, we pick independently a bijection $\pi(\q;\cdot)$ with a 
uniform law over the set of bijections satisfying the reversibility and no U--turn 
constraints. For $d\geq 2$, there are $(2d-1)!!$ such bijections, which is also the 
number of pairings of $2d$ edges. We denote by $\bbM$ the corresponding probability 
law over the mirrors $\{\pi(\q,\cdot):q\in \Lambda_N\}$, and by $\langle \cdot \rangle$ 
the average with respect to this law.

We define, using the notation $x=(\q,\p)$ for a generic point of $\cM_N$
\begin{equation}
T(x,x';\pi)=
\begin{cases}
1 & \text{if } F(x,\pi)=x',\\
0 & \text{otherwise}.
\end{cases}
\label{defT}
\end{equation}

Since $F$ is a bijection, for each $x$, there is exactly one $x'$ such that 
$T(x,x';\pi)=1$, and we have
\begin{equation}
\sum_{x'\in\cM_N\cup O_N}T(x,x';\pi)=1,
\quad \forall x,\ \forall \pi.
\label{sumrule}
\end{equation}

The slab $\Lambda_N$ itself may be viewed as a single deterministic scatterer: for 
a fixed configuration of mirrors, each incoming point $x\in I_N$ is mapped to a 
unique outgoing point in $O_N$. For any $x\in I_N$, define the exit time
\[
n_x:=\inf \{n\ge0: F^n(x;\pi)\in O_N\},
\]
and the induced map
\begin{equation}
F_N(x;\pi)=F^{n_x}(x;\pi).
\label{Fboxm}
\end{equation}
Define
\begin{equation}
T_N(x,x';\pi)=
\begin{cases}
1 & \text{if } F_N(x,\pi)=x',\\
0 & \text{otherwise}.
\end{cases}
\label{defT1}
\end{equation}
Since $F$ is a bijection between $ \cM_N$ and $(\cM_N\backslash I_N)\cup O_N$ and 
$\Lambda_N$ is finite, for any $x\in I_N$, $n_x$ is finite and not larger than 
$2d|\Lambda_N|$. There can be no loop containing a point in $I_N$ because no point 
$x\in \cM_N$ is mapped to $I_N$.

The single-point crossing probability is
\begin{equation}
p_N(x,x'):=\left \langle T_N(x,x';\pi)\right\rangle.
\label{crossing_prob}
\end{equation}
Our main object of study is the probability that a particle entering the volume 
$\Lambda_N$ on the left side at point ${\bf 1}:=(1,\ldots,1)$ crosses the volume 
and exits at some site $(N+1,q_2,\ldots,q_d)$:
\[
C_N:=\sum_{x'\in O^+_N}\left \langle T_N(x,x';\pi)\right \rangle,
\]
with $x=(\mathbf{1},\e_1)$. We also introduce the conductivity at scale $N$:
\[
\kappa_N=\frac{N C_N}{1-C_N}.
\]
We will argue that $\kappa_N\to\kappa_\infty$ as $N\to\infty$ with $0<\kappa_\infty<\infty$.


\section{Markovian Benchmark: Non--Backtracking Random Walk}

We now consider the average of the transition rule with respect to the disorder:
\begin{equation}
p(x,x')=\langle T(x,x';\pi)\rangle
=\frac 1 {2d-1}\delta_{\q',\q+\p'}(1-\delta_{\p,-\p'}),
\label{NBWkernel}
\end{equation}
which is a transition probability for a Markov chain $(X_n)_{n\in\bbN}$ on 
$\cM_N$. This is the transition kernel of a kinematic non--backtracking random walk.

Let $\widehat C_N$ be the corresponding crossing probability:
\[
\widehat C_N=\bbP[X_{t_x}\in O^+_N\mid X_0=x],
\]
for any $x\in I^+_N$, where $t_x=\inf\{n: X_n\in O_N\}$ is the exit time.

It is possible to compute $\widehat C_N$ explicitly and show
\begin{equation}
\widehat C_N=\frac{\widehat\kappa_0}{\widehat\kappa_0+N}, 
\qquad \widehat\kappa_0=\frac{d}{d-1}.
\end{equation}
To see this, define $\widehat c_n=\widehat C_{2^n}$ and
\[
\widehat \kappa_n=\frac{2^n\widehat c_n}{1-\widehat c_n}.
\]
A slab of length $N=2^{n+1}$ can be divided into two equal parts of length $2^n$, 
say $\Lambda_1$ and $\Lambda_2$. A path starting at $x\in I^+_{2^{n+1}}$ can exit in 
$O^+_{2^{n+1}}$ by crossing first $\Lambda_1$ and then $\Lambda_2$, or by revisiting 
the interface between them.

Independence yields
\begin{equation}
\widehat c_{n+1} = \sum_{l=0}^\infty (\widehat c_n)^2(1-\widehat c_n)^{2l}
=\frac{\widehat c_n}{2-\widehat c_n}.
\label{iterateM1}
\end{equation}
It is easy to check that $\widehat \kappa_{n+1}=\widehat \kappa_n$, and thus 
$\widehat \kappa_n=\widehat \kappa_0$ for all $n\geq 0$, and
\[
\widehat c_n=\frac{\widehat \kappa_0}{2^n+\widehat \kappa_0}.
\]
Computing $\widehat c_0$ gives
\[
\widehat c_0=\frac 1 {2d-1}+\frac{2d-2}{2d-1}\frac{1}{2}=\frac{d}{2d-1},
\]
and hence $\widehat \kappa_0=d/(d-1)$.


\section{Multiscale Decomposition and Recursion}
\begin{figure}[t]
\centering
\begin{tikzpicture}[scale=0.50,>=Latex]

\Slab{(0,0)}{2.2}{blue}
\Slab{(2.3,0)}{2.2}{blue}

\draw[dashed,black!70] (2.3,0) -- (2.3,0.9)
                       (3.2,0.45) -- (3.2,1.35);

\draw[decorate,decoration={brace,mirror,amplitude=5pt}]
      (0,-0.5) -- node[below=4pt] {$2^{n}$} (2.2,-0.5);
\draw[decorate,decoration={brace,mirror,amplitude=5pt}]
      (2.3,-0.5) -- node[below=4pt] {$2^{n}$} (4.5,-0.5);

\draw[very thick,->] (5.0,1.1) -- (6.6,1.1) node[midway,above]{stack};

\Slab{(6.9,0)}{4.5}{teal}
\draw[decorate,decoration={brace,mirror,amplitude=5pt}]
      (6.9,-0.5) -- node[below=4pt] {$2^{\,n+1}$} (11.4,-0.5);

\Slab{(0,-4.2)}{4.5}{teal}
\Slab{(4.7,-4.2)}{4.5}{teal}

\draw[dashed,black!70] (4.7,-4.2) -- (4.7,-3.3)
                       (5.6,-3.75) -- (5.6,-2.85);

\draw[decorate,decoration={brace,mirror,amplitude=5pt}]
      (0,-4.7) -- node[below=4pt] {$2^{\,n+1}$} (4.5,-4.7);
\draw[decorate,decoration={brace,mirror,amplitude=5pt}]
      (4.7,-4.7) -- node[below=4pt] {$2^{\,n+1}$} (9.2,-4.7);

\draw[very thick,->] (9.8,-3.1) -- (11.6,-3.1) node[midway,above]{stack};

\Slab{(12.0,-4.2)}{9.0}{orange}
\draw[decorate,decoration={brace,mirror,amplitude=5pt}]
      (12.0,-4.7) -- node[below=4pt] {$2^{\,n+2}$} (21.0,-4.7);


\end{tikzpicture}
\caption{Multiscale composition in 3D: $2^n+2^n\!\to\!2^{n+1}$, then $2^{n+1}+2^{n+1}\!\to\!2^{n+2}$.}
\end{figure}

We now turn to the mirrors model, where the dynamics is non--Markovian.  
Define $c_n=C_{2^n}$ and 
\[
\kappa_n=\frac{2^n c_n}{1-c_n}.
\]

The point of introducing $\kappa_n$ is that the normal conduction scaling
$c_n\sim \kappa/2^n$ is equivalent to the existence of a finite nonzero limit of
$\kappa_n$. Our recursive analysis is organized precisely around this quantity.
We regard $c_n$ as already being of order $2^{-n}$, and we use the concatenation
of two slabs of size $2^n$ to compute the correction to the effective conductivity
when passing to scale $2^{n+1}$. Thus the multiscale step is naturally formulated
as a recursion for $\kappa_n$, rather than only for $c_n$.
As above, we split a slab of length $2^{n+1}$ into two equal parts of length $2^n$, 
denoted $\Lambda_1$ (left) and $\Lambda_2$ (right). For each part we define the 
incoming phase-space sets $\mathcal I_i^{\pm}$ and outgoing phase-space sets 
$\mathcal O_i^{\pm}$ in a way strictly analogous to the definitions in Section~2.
The two half-slabs are separated by the interface between $\mathcal O_1^+$ and 
$\mathcal O_2^-$. Since both describe the same geometric interface with opposite 
orientations, they are naturally in one-to-one correspondence.

For a slab of width $2^n$, we let $t_n=T_{2^n}$ denote the indicator of a 
\emph{boundary transfer} at scale $2^n$. More precisely, for an incoming boundary 
state $a\in \mathcal I_n:=\mathcal I_n^+\cup \mathcal I_n^-$ and an outgoing 
boundary state $b\in \mathcal O_n:=\mathcal O_n^+\cup \mathcal O_n^-$, we define
\[
t_n(a,b;\pi)=
\begin{cases}
1 & \text{if }F_{2^n}(a;\pi)=b,\\
0 & \text{otherwise.}
\end{cases}
\]
Thus $t_n(a,b;\pi)$ records the full transfer through the slab of width $2^n$; 
it is \emph{not} restricted to crossings from one side to the other, and it also includes transfers that exit on 
the same side as the entrance. In fact if $a\in \cI^+_n$ and $b\in\cO^-_n$ or $a\in \cI^-_n$ and $b\in\cO^+_n$, we call this type of transfer a {\it return}. Its expectation is
\begin{equation}
c_n(a,b):=\langle t_n(a,b;\pi)\rangle.
\label{proban}
\end{equation}
The crossing probability $c_n=C_{2^n}$ is the special case obtained by summing 
$c_n(a,b)$ over right exits when the initial state $a=(\un 1,\e_1)$ belongs to 
$\mathcal I_n^+$.

A path crossing the full slab of length $2^{n+1}$ may alternate several times 
between the two half-slabs before its final exit on the right. We consider an integer $l$ that
counts the total number of transfers on each side of the interface in the concatenation scheme. 
Equivalently, after the first crossing of the left half-slab and reaching the interface, the trajectory makes 
$l-1$ further interface returns on each side of it, before exiting the full slab on the right. Thus 
$l=1$ corresponds to a direct crossing of the two half-slabs, while $l=2$ 
corresponds to one interface return before the final exit.

To write this decomposition explicitly, let the initial state be
\[
x_1=x=(\un 1,\e_1)\in \mathcal I_1^+,
\]
and let $x_{l+1}=x'\in \mathcal O_2^+$ be the final exit state on the right. For 
$j=1,\dots,l$, let $y_j\in \mathcal O_1^+$ denote the successive transfer states 
from the left half-slab to the interface, and for $j=2,\dots,l$, let 
$x_j\in \mathcal O_2^-$ denote the successive transfer states from the right 
half-slab back to the interface. Then, for a fixed environment $\pi=(\pi_1,\pi_2)$, 
a path contributing to the crossing of the full slab with exactly $l-1$ interface 
returns is represented by the ordered sequence of transfers
\[
x_1 \xrightarrow{\,\Lambda_1\,} y_1 \xrightarrow{\,\Lambda_2\,} x_2
\xrightarrow{\,\Lambda_1\,} y_2 \xrightarrow{\,\Lambda_2\,} \cdots
\xrightarrow{\,\Lambda_1\,} y_l \xrightarrow{\,\Lambda_2\,} x_{l+1}.
\]
Therefore
\[
t_{n+1}(x,x';\pi)
= \sum_{l\geq 1}
\sum_{\substack{y_1,\dots,y_l\in \mathcal O_1^+\\ x_2,\dots,x_l\in \mathcal O_2^-}}
\prod_{i=1}^{l}t_n(x_i,y_i;\pi_1)
\prod_{i=1}^l t_n(y_i,x_{i+1};\pi_2).
\]
This is the expanded version of the concatenation rule: the first product collects 
the transfers through the left half-slab, the second product collects the transfers 
through the right half-slab, and all intermediate interface states are summed over.

Using independence of $\pi_1$ and $\pi_2$ and taking expectation yields
\begin{equation}
c_{n+1}
=\sum_{l\geq 1}
\sum_{x_{l+1}\in \mathcal O^+_2}
\sum_{\substack{y_1,\dots,y_l\in \mathcal O_1^+\\ x_2,\dots,x_l\in \mathcal O_2^-}}
\left\langle \prod_{i=1}^{l}t_n(x_i,y_i;\pi)\right\rangle 
\left\langle\prod_{i=1}^l t_n(y_i,x_{i+1};\pi)\right\rangle,
\label{recursion}
\end{equation}
with $x_1=(\un 1,\e_1)$. By translation invariance in the 
directions orthogonal to ${\mathbf e}_1$, the term corresponding to $l=1$ is 
$(c_n)^2$.

We introduce
\begin{equation}
\eta_n(l)=\frac{1}{c_n^2(1-c_n)^{2(l-1)}}
\sum_{x_{l+1}\in \mathcal O^+_2}
\sum_{\substack{y_1,\dots,y_l\in \mathcal O_1^+\\ x_2,\dots,x_l\in \mathcal O_2^-}}
\left\langle\prod_{i=1}^l t_n(x_i,y_i;\pi)\right\rangle
\left\langle\prod_{i=1}^l t_n(y_i,x_{i+1};\pi)\right\rangle.
\label{etadef}
\end{equation}
This measures how correlations at scale $2^n$ deviate from independent concatenation 
for trajectories involving $l$ transfer blocks, equivalently $l-1$ further returns to the interface after the first arrival.

We write
\[
c_{n+1}=\sum_{l=1}^\infty c_n^2(1-c_n)^{2(l-1)}\eta_n(l).
\]
Defining
\[
S_n=\sum_{l=1}^\infty (1-c_n)^{2(l-1)}\eta_n(l),
\]
with $\eta_n(1)=1$, we have
\[
c_{n+1}=c_n^2 S_n.
\]

We also observe
\[
S_n-(1-c_n)^2S_n
=1+\sum_{l=1}^\infty (1-c_n)^{2l}(\eta_n(l+1)-\eta_n(l)),
\]
hence
\[
S_n=\frac{1}{1-(1-c_n)^2}(1+\Delta_n),
\]
with
\begin{equation}
\Delta_n=\sum_{l=1}^\infty (1-c_n)^{2l}(\eta_n(l+1)-\eta_n(l)).
\label{Deltan}
\end{equation}
The analog of \eqref{iterateM1} becomes
\[
c_{n+1}=\frac{c^2_n(1+\Delta_n)}{1-(1-c_n)^2}
=\frac{c_n}{2-c_n}(1+\Delta_n).
\]

We are interested in the induced recursion for $\kappa_n$:
\begin{equation}
\kappa_{n+1}=\frac{2^{n+1} c_{n+1}}{1-c_{n+1}}
=\frac{2^n c_n(1+\Delta_n)}{1-c_n-\frac{c_n}{2}\Delta_n}.
\label{recursionkappa2}
\end{equation}

Equation \eqref{recursionkappa2} is the basic recursive relation of the paper:
assuming normal conduction scaling at scale $2^n$, it expresses the correction to the
conductivity constant at scale $2^{n+1}$ through the correlation term $\Delta_n$.
More precisely, the normal conduction hypothesis $c_n\asymp 2^{-n}$ is not used to derive
the decomposition itself, but only to evaluate the size of the correction terms
appearing in $\Delta_n$ and hence in the recursion for $\kappa_n$.


\section{Closure Assumption and Second-Order Terms}

In the definition (\ref{etadef}) of $\eta_n(l)$ and hence $\Delta_n$, each term is the product of two 
factors, each representing the joint probability of $l$ boundary transfers in a 
box of size $2^n$, with prescribed incoming and outgoing states. In other words, the 
crossing probability at scale $n+1$ depends on correlations between several transfer 
events at scale $n$ occurring inside the same half-slab.

The basic reason these transfer events are not independent is that a slab of width 
$2^n$ defines, for every mirror configuration, a bijection between incoming and outgoing 
boundary states. This bijectivity, together with reversibility, produces two distinct 
effects on multipoint functions. First, there are \emph{hard constraints}: some pairs of 
prescribed transfer events are impossible, while others are forced. Second, once these 
hard constraints are removed, one expects the residual correlation between admissible 
transfer events to become weak as the scale grows. The closure assumption used below 
formalizes precisely this asymptotic factorization. See Appendix~B for the exhaustive 
case-by-case structural formulas.

To make this explicit without overloading the notation of the recursion, we denote in 
this discussion by $a_1,a_2$ generic incoming boundary states and by $b_1,b_2$ generic 
outgoing boundary states of a single slab of width $2^n$. We first associate to each 
boundary transfer $(a,b)$ two labels.

The first one records whether the transfer is a crossing or a return:
\[
\chi(a,b)=
\begin{cases}
c,& \text{if } a\in I^+,\, b\in O^+ \text{ or } a\in I^-,\, b\in O^-,\\[1mm]
r,& \text{if } a\in I^+,\, b\in O^- \text{ or } a\in I^-,\, b\in O^+.
\end{cases}
\]
Thus $\chi(a,b)=c$ means that the transfer exits on the opposite side of the slab, while 
$\chi(a,b)=r$ means that it returns to the same side.

The second label records on which side a return occurs. We define
\[
\sigma(a,b)=
\begin{cases}
L,& \text{if } a\in I^+,\, b\in O^-,\\[1mm]
R,& \text{if } a\in I^-,\, b\in O^+,\\[1mm]
\mathrm{cross},& \text{if } \chi(a,b)=c.
\end{cases}
\]
The role of $\sigma$ is to distinguish the two geometrically different return--return 
configurations: one return on each side and two returns on the same side.

We also define the forced-coincidence indicator
\[
\mathbf 1_{\mathrm{for}}(a_1,b_1,a_2,b_2)
:=
\delta_{a_1a_2}\delta_{b_1b_2}
+\delta_{a_1Rb_2}\delta_{a_2Rb_1},
\]
which encodes the two deterministic ways in which two prescribed transfers can actually 
coincide: either they are literally the same transfer, or they are paired by time reversal.

Next we define the admissibility indicator. The admissible sector depends on the transfer 
types and on the side label of return transfers. Writing
\[
\chi_i:=\chi(a_i,b_i),\qquad \sigma_i:=\sigma(a_i,b_i),\qquad i=1,2,
\]
we set
\begin{align}
\mathbf 1_{\mathrm{adm}}(a_1,b_1,a_2,b_2)
&= \mathbf 1_{\{\chi_1=\chi_2=r,\ \sigma_1\neq \sigma_2\}}
\notag\\
&\quad + \mathbf 1_{\{\chi_1=\chi_2=r,\ \sigma_1=\sigma_2\}}
(1-\delta_{a_1a_2})(1-\delta_{b_1b_2})(1-\delta_{a_1Rb_2})(1-\delta_{a_2Rb_1})
\notag\\
&\quad + \mathbf 1_{\{(\chi_1,\chi_2)=(c,r)\ \text{or}\ (r,c)\}}
(1-\delta_{Rb_1a_2})(1-\delta_{b_1b_2})(1-\delta_{a_1a_2})(1-\delta_{a_1Rb_2})
\notag\\
&\quad + \mathbf 1_{\{\chi_1=\chi_2=c\}}
(1-\delta_{b_1b_2})(1-\delta_{a_1a_2})(1-\delta_{Rb_2a_1})(1-\delta_{Rb_1a_2}).
\label{admindicator}
\end{align}
This is simply a compact way of packaging the five geometric sectors listed explicitly in 
Appendix~B. In particular:
\begin{itemize}
\item $\chi=\chi'=r$ and $\sigma\neq \sigma'$ corresponds to $A^1_{rr}$;
\item $\chi=\chi'=r$ and $\sigma=\sigma'$ corresponds to $A^2_{rr}$;
\item $(\chi,\chi')=(c,r)$ corresponds to $A_{cr}$;
\item $(\chi,\chi')=(r,c)$ corresponds to $A_{rc}$;
\item $\chi=\chi'=c$ corresponds to $A_{cc}$.
\end{itemize}

With this notation, the closure assumption takes the following compact form:
\begin{equation}
\begin{aligned}
\langle t_n(a_1,b_1)t_n(a_2,b_2)\rangle
&=
\mathbf 1_{\mathrm{adm}}(a_1,b_1,a_2,b_2)\,
c_n(a_1,b_1)c_n(a_2,b_2)
\\
&\qquad\cdot \bigl(1+r_n(a_1,b_1,a_2,b_2)\bigr)
+\mathbf 1_{\mathrm{for}}(a_1,b_1,a_2,b_2)\,c_n(a_1,b_1).
\end{aligned}
\label{masterclosure}
\end{equation}
where the remainder satisfies the uniform bound
\[
|r_n(a_1,b_1,a_2,b_2)|\le h(n),
\qquad h(n)\xrightarrow[n\to\infty]{}0.
\]
Thus $h(n)$ is a uniform bound on the relative closure error. It is not an additional 
parameter to be fitted. Rather, it is a bookkeeping function that quantifies the uniform 
error of the closure approximation. The analysis below only uses the fact that 
$h(n)\to 0$ as $n\to\infty$. In the leading-order computation of $\alpha$, one neglects 
the terms controlled by $h(n)$, that is, effectively sets them to zero at first order. 

In words, the two-point function is the product of one-point functions on the admissible sector, up to a relative error bounded by $h(n)$, while the forced sector corresponds to deterministic coincidences induced by bijectivity and time reversal.

The compact formula \eqref{masterclosure} summarizes all the structural cases used later.
For the specific $l=2$ contribution to the recursion, the relevant two-point functions
belong first to the mixed sectors $A_{cr}$ and $A_{rc}$. 

We now return to the recursion variables and analyze the case $l=2$, which gives the 
dominant correction. In the expression (\ref{etadef}) for $\eta_n(l)$ we insert
\begin{equation}
t_n(x,y;\pi)=c_n(x,y)+\delta t_n(x,y),
\label{insertion}
\end{equation}
for each factor in the products. This \emph{defines} $\delta t_n$. The term with no 
$\delta t$ gives 1, so it is natural to define $\delta_n(l)=\eta_n(l)-1$.

Assuming absolute convergence of the series defining $\Delta_n$ and using 
$\eta_n(1)=1$, we find
\begin{eqnarray}
\Delta_n
&=&\sum_{l=1}^\infty \left((1-c_n)^{2l}-(1-c_n)^{2(l+1)}\right)\delta_n(l+1)\nonumber\\
&=& c_n(2-c_n)\sum_{l=1}^\infty (1-c_n)^{2l}\delta_n(l+1).	
\end{eqnarray}

For $l=2$, the definition (\ref{etadef}) of $\eta_n(2)$ reads
\[
\eta_n(2)=\frac{1}{c_n^2(1-c_n)^2}
\sum_{y_1\in O_1^+}\sum_{x_2\in O_2^-}\sum_{y_2\in O_1^+}\sum_{x_3\in O_2^+}
\left\langle t_n(x_1,y_1)t_n(x_2,y_2)\right\rangle
\left\langle t_n(y_1,x_2)t_n(y_2,x_3)\right\rangle,
\]
with $x_1=(\un 1,\e_1)$ fixed. 
We now insert
\[
t_n(u,v;\pi)=c_n(u,v)+\delta t_n(u,v),
\qquad \langle \delta t_n(u,v;\pi)\rangle=0,
\]
into the two expectations above. Expanding the product gives four types of terms. The term
with no factor $\delta t_n$ is exactly
\[
\frac{1}{c_n^2(1-c_n)^2}
\sum_{y_1,x_2,y_2,x_3}
c_n(x_1,y_1)c_n(x_2,y_2)c_n(y_1,x_2)c_n(y_2,x_3)=1,
\]
which is cancelled by the subtraction defining $\delta_n(2)=\eta_n(2)-1$. The terms containing exactly
one centered factor $\delta t_n$ vanish because their expectation is zero. Thus only the
terms containing two centered factors survive. Grouping them according to whether the two
centered factors lie in both expectations or in exactly one of them, we obtain
\[
\delta_n(2)=R_1+R_2,
\]
where
\begin{align}
c_n^2(1-c_n)^2R_1
&=\sum_{y_1\in O_1^+}\sum_{x_2\in O_2^-}\sum_{y_2\in O_1^+}\sum_{x_3\in O_2^+}
\langle\delta t_n(x_1,y_1)\delta t_n(x_2,y_2)\rangle
\langle\delta t_n(y_1,x_2)\delta t_n(y_2,x_3)\rangle,
\label{R1}
\\
c_n^2(1-c_n)^2R_2
&=2\sum_{y_1\in O_1^+}\sum_{x_2\in O_2^-}\sum_{y_2\in O_1^+}\sum_{x_3\in O_2^+}
\langle\delta t_n(x_1,y_1)\delta t_n(x_2,y_2)\rangle 
c_n(y_1,x_2)c_n(y_2,x_3).
\label{R2}
\end{align}
Here $x_1=(\un 1,\e_1)$ is fixed.

We now state the closure assumption in the case $l=2$ and in the special case appearing 
in (\ref{R1}) and (\ref{R2}). In these expressions, the first two-point function
\[
\langle \delta t_n(x_1,y_1)\delta t_n(x_2,y_2)\rangle
\]
belongs to the crossing--return sector $A_{cr}$: indeed, $x_1\in I_1^+$ and 
$y_1\in O_1^+$ so that $x_1\to y_1$ is a crossing transfer through $\Lambda_1$, whereas 
$x_2\in O_2^-$ and $y_2\in O_1^+$ correspond, after identification of the interface, to a 
return transfer through $\Lambda_1$. Similarly, the second two-point function
\[
\langle \delta t_n(y_1,x_2)\delta t_n(y_2,x_3)\rangle
\]
belongs to the return--crossing sector $A_{rc}$: the transfer $y_1\to x_2$ is a return 
through $\Lambda_2$, while $y_2\to x_3$ is a crossing through $\Lambda_2$. Therefore, 
in the present situation, \eqref{masterclosure} reduces to the $A_{cr}$ and $A_{rc}$ 
forms.

More explicitly, for the first expectation one has
\begin{align}
\langle t_n(x_1,y_1)t_n(x_2,y_2)\rangle
&=(1-\delta_{Ry_1x_2})(1-\delta_{y_1y_2})(1-\delta_{x_1x_2})(1-\delta_{x_1Ry_2})
\notag\\
&\qquad\cdot c_n(x_1,y_1)c_n(x_2,y_2)\bigl(1+O(h(n))\bigr),
\label{closurecr}
\end{align}
while for the second expectation one has
\begin{align}
\langle t_n(y_1,x_2)t_n(y_2,x_3)\rangle
&=(1-\delta_{x_2y_2})(1-\delta_{x_2x_3})(1-\delta_{y_1y_2})(1-\delta_{y_1Rx_3})
\notag\\
&\qquad\cdot c_n(y_1,x_2)c_n(y_2,x_3)\bigl(1+O(h(n))\bigr).
\label{closurerc}
\end{align}
Since in (\ref{R1}) and (\ref{R2}) the initial state $x_1$ lies on the left boundary, the
factors $(1-\delta_{x_1x_2})$ and $(1-\delta_{x_1Ry_2})$ are identically equal to $1$.
Hence the first closure relation reduces to
\begin{equation}
\langle t_n(x_1,y_1)t_n(x_2,y_2)\rangle
=(1-\delta_{Ry_1x_2})(1-\delta_{y_1y_2})\,c_n(x_1,y_1)c_n(x_2,y_2)\bigl(1+O(h(n))\bigr).
\label{closure}
\end{equation}
where $h(n)\to 0$ as $n\to \infty$. 
A similar simplification holds for the second two-point function in \eqref{R1}: since
$x_3\in O_2^+$ is a right exit state, the factor $(1-\delta_{x_2x_3})$ is identically
equal to $1$, while the remaining Kronecker factors encode the residual deterministic
constraints. In the subsequent estimates, these constraints are either used explicitly or
eliminated by the summation identities described below.
The important point for the present derivation is 
that the only non-factorized contributions are those imposed by the deterministic 
constraints, while the remaining admissible sector factorizes up to a uniform error 
controlled by $h(n)$.

Using (\ref{closure}), one can show (see Appendix~C) that
\begin{equation}
|R_2|\leq c_n\,O(h(n)),
\label{R2bound}
\end{equation}
consistent with numerical simulations, which indicate a decay of order $1/(2^n)^2$.
At this point the normal conduction scaling enters explicitly: under the inductive hypothesis
$c_n\asymp 2^{-n}$, the bound \eqref{R2bound} becomes
\[
|R_2|=2^{-n}O(h(n)).
\]
Thus as we will see below, within the \(l=2\) contribution to \(\Delta_n\), the term \(R_2\) is smaller than the term \(R_{1}\) by an extra factor \(c_n\), and is therefore negligible as soon as \(h(n)\to0\).
The smallness of $R_2$ has a simple structural origin. The first two-point function in
\eqref{R2} belongs to the mixed sector $A_{cr}$, which a priori is not smaller than the
dominant contribution. However, the sum over the second outgoing variable may be converted,
using the deterministic transfer identity : for any $x\in \cI_1$,
\[
\sum_{y\in O_1^+} t_n(x,y)+\sum_{y\in O_1^-} t_n(x,y)=1,
\]
and its averaged version for $c_n$, into a complementary contribution involving transfers
to $O_1^-$. In this way the average of the mixed crossing--return structure is turned into a
crossing--crossing one. The latter is smaller because it is approximately (under the closure assumption) equal to the product of two crossing probabilities. Thus the negligibility of $R_2$ is a
consequence of the deterministic bijective dynamics together with the fact that, after this
rearrangement, one is left with a sector involving an additional crossing constraint.

We now turn to the dominant term $R_1$. Writing $R_1=R_{11}+R_{12}$ and splitting 
the sum in (\ref{R1}) according to $y_2=y_1$ or $y_2\neq y_1$, the terms with 
$y_2\neq y_1$ can be shown (again using closure) to be of order $(h(n))^2$ and thus 
negligible. 
The dominant contribution comes from the diagonal sector $y_2=y_1$. Since $R_1$ is
defined via centered variables $\delta t_n=t_n-c_n$, we must evaluate
\[
\langle \delta t_n(x_1,y_1)\,\delta t_n(x_2,y_1)\rangle
\qquad\text{and}\qquad
\langle \delta t_n(y_1,x_2)\,\delta t_n(y_1,x_3)\rangle .
\]

For the first factor, by bijectivity of the slab transfer map, two distinct incoming
states cannot be mapped to the same outgoing state. Since $x_1\in\cI_1^+$ and
$x_2\in\cO_2^-$ lie on different boundaries, one has $x_1\neq x_2$, hence
$t_n(x_1,y_1)t_n(x_2,y_1)=0$ for every mirror configuration. Therefore
\[
\langle \delta t_n(x_1,y_1)\,\delta t_n(x_2,y_1)\rangle
= -c_n(x_1,y_1)\,c_n(x_2,y_1).
\]

For the second factor, since $x_2\in\cO_2^-$ and $x_3\in\cO_2^+$, one has
$x_2\neq x_3$. Again by bijectivity, a given incoming state $y_1$ cannot be mapped
simultaneously to two distinct outgoing states, so
$t_n(y_1,x_2)t_n(y_1,x_3)=0$. Hence
\[
\langle \delta t_n(y_1,x_2)\,\delta t_n(y_1,x_3)\rangle
= -c_n(y_1,x_2)\,c_n(y_1,x_3).
\]

Multiplying the two contributions and using time-reversal invariance together with the interface identification,
$c_n(x_2,y_1)=c_n(Ry_1,Rx_2)=c_n(y_1,x_2)$, we obtain
\[
c_n^2(1-c_n)^2R_{11}
= \sum_{y_1\in\cO_1^+}\sum_{x_2\in\cO_2^-}\sum_{x_3\in\cO_2^+}
  c_n(x_1,y_1)\,\bigl(c_n(y_1,x_2)\bigr)^2\,c_n(y_1,x_3).
\]
By translation invariance 
in $y_1$, we have $\sum_{x_3\in O_2^+} c_n(y_1,x_3)=c_n$, and hence
\[
c_n(1-c_n)^2 R_{11}
=\sum_{y_1\in O_1^+}\sum_{x_2\in O_2^-}
c_n(x_1,y_1)\bigl(c_n(y_1,x_2)\bigr)^2.
\]
Using again translation invariance in $y_1$, this gives
\[
(1-c_n)^2 R_{11}
=\sum_{x_2\in O_2^-}\bigl(c_n(y_1,x_2)\bigr)^2
\leq \Bigl(\sum_{x_2\in O_2^-}c_n(y_1,x_2)\Bigr)^2
=(1-c_n)^2.
\]
Thus $R_{11}\leq 1$. Numerical computation shows rapid convergence of 
$(1-c_n)^2R_{11}$ to a constant independent of $n$:
\[
(1-c_n)^2R_{11}\to 0.018704\pm 10^{-5}.
\]
The contribution to $\Delta_n$ of the term $l=2$ reads
\[
\Delta_n(2)= c_n(2-c_n)\bigl((1-c_n)^2R_{11}+O(h(n))\bigr).
\]
This is the second place where the normal conduction scaling is used. Indeed, under the
inductive assumption $c_n\asymp 2^{-n}$ and since $(1-c_n)^2R_{11}$ converges to a
finite nonzero constant, the whole contribution $\Delta_n(2)$ is of order $2^{-n}$.
Equivalently, using $c_n=\kappa_n/(2^n+\kappa_n)$ with $\kappa_n=O(1)$, one gets
\[
\Delta_n(2)=\frac{\kappa_n}{2^n}\bigl(\alpha+o(1)\bigr),
\]
which is precisely the scale of the correction entering the recursion for
$\kappa_{n+1}$.
\begin{figure}[h]
\centering
\begin{tikzpicture}[scale=0.3,every node/.style={font=\small}]
  \draw[step=1cm,lightgray,very thin] (0,0) grid (14,14);

  \draw[thick] (0,0) rectangle (7,14);
  \draw[thick] (7,0) rectangle (14,14);

  \draw[blue,thick]
    (0,1) -- (0,2) -- (1,2) -- (2,2) -- (2,3) -- (3,3)
           -- (4,3) -- (4,4) -- (5,4) -- (6,4) -- (7,4) -- (7,5);

  \draw[blue,thick]
    (7,8) -- (6,8) -- (5,8) -- (5,7) -- (5,6) -- (6,6)
          -- (6,5) -- (7,5);

  \draw[green,thick]
    (7,5) -- (8,5) -- (8,6) -- (9,6) -- (9,7) -- (8,7)
          -- (8,8) -- (7,8);

  \draw[green,thick]
    (7,5) -- (8,5) -- (9,5) -- (10,5) -- (11,5) -- (12,5)
          -- (12,4) -- (12,3) -- (13,3) -- (14,3) -- (14,2);

  \node[below left] at (0,1)   {$x_1$};
  \node[below ] at (7,5)  {$y_1=y_2$};
  \node[above right] at (7,8)  {$x_2$};
  \node[below right] at (14,2) {$x_3$};

\end{tikzpicture}
\caption{Trajectories contributing to $R_{11}$ for $l=2$.}
\end{figure}

\begin{figure}[htbp]
\centering
\includegraphics[width=.50\textwidth]{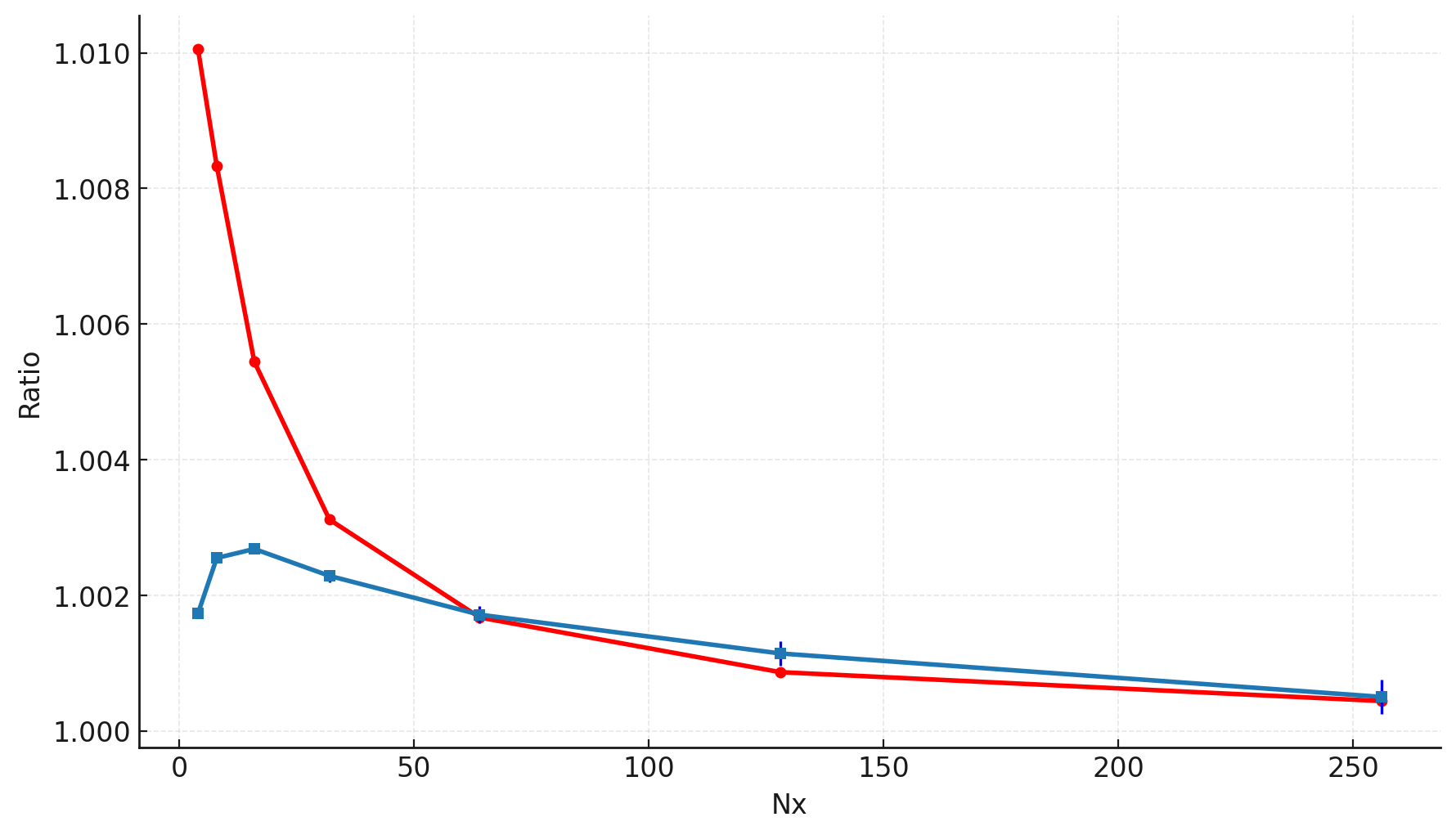}
\caption{The measured ratio $1+\Delta_n=c_{n+1}(2-c_n)/c_n$ is plotted in blue with 
a 95\% confidence interval, while $1+c_n(2-c_n)(1-c_n)^2R_{11}$ is plotted in red.  $c_n$ is estimated numerically by sampling random mirror environments, launching a particle from a fixed boundary point, and recording the frequency with which it exits on the opposite side.  As explained in the main text $R_{11}$ is obtained via the sum $\sum_{x_2\in O_2^-}(c_n(y_1,x_2))^2$.  The ``return'' probabilities $c_n(y_1,x_2)$ are estimated numerically, and the sum of their square converges robustly (almost independently  of the system size). As explained in the appendix, it is possible to compute systematically lower bounds on those return probabilities.  It is also likely that one can write an expansion on paths for those objects.
The contribution of $R_2$ lowers the red values significantly only for the two first 
values of $N$, the contribution of $R_{12}$ increases those values in a way that is 
almost not visible at this scale. }
\label{fig:ratiosR11}
\end{figure}

At this stage, the leading factor $c_n^2$ in the recursion for $c_{n+1}$ simply
reflects the fact that crossing the large slab requires crossing both half-slabs.
If $c_n$ is of order $2^{-n}$, then $c_{n+1}$ is again of order
$2^{-(n+1)}$, and the nontrivial information lies in the prefactor. This is why
we rewrite the recursion in terms of $\kappa_n$: the scale doubling preserves the
leading conductive order, while the correlation term $\Delta_n$ produces the
renormalization of the conductivity constant.
Using \eqref{recursionkappa2} and $c_n=\kappa_n/(2^n+\kappa_n)$, we obtain
\[
\kappa_{n+1}=\kappa_n\left(1+\frac{\kappa_n}{2^n} (\alpha+o(1))\right),
\qquad \alpha=0.0374,
\]
where $\alpha$ is twice the limiting value of $(1-c_n)^2R_{11}$. Iterating this relation 
(neglecting $o(1)$) from $\kappa_n=1.5397\pm 3 \times 10^{-4}$ at $n=8$ yields
\[
\kappa_\infty=1.5403.
\]
The comparison between the measured ratio 
\[
1+\Delta_n=\frac{c_{n+1}(2-c_n)}{c_n}
\]
and the result of the above computation is given in Figure~\ref{fig:ratiosR11}.

\begin{figure}[h!]
\begin{center}
\includegraphics[width=.40\textwidth]{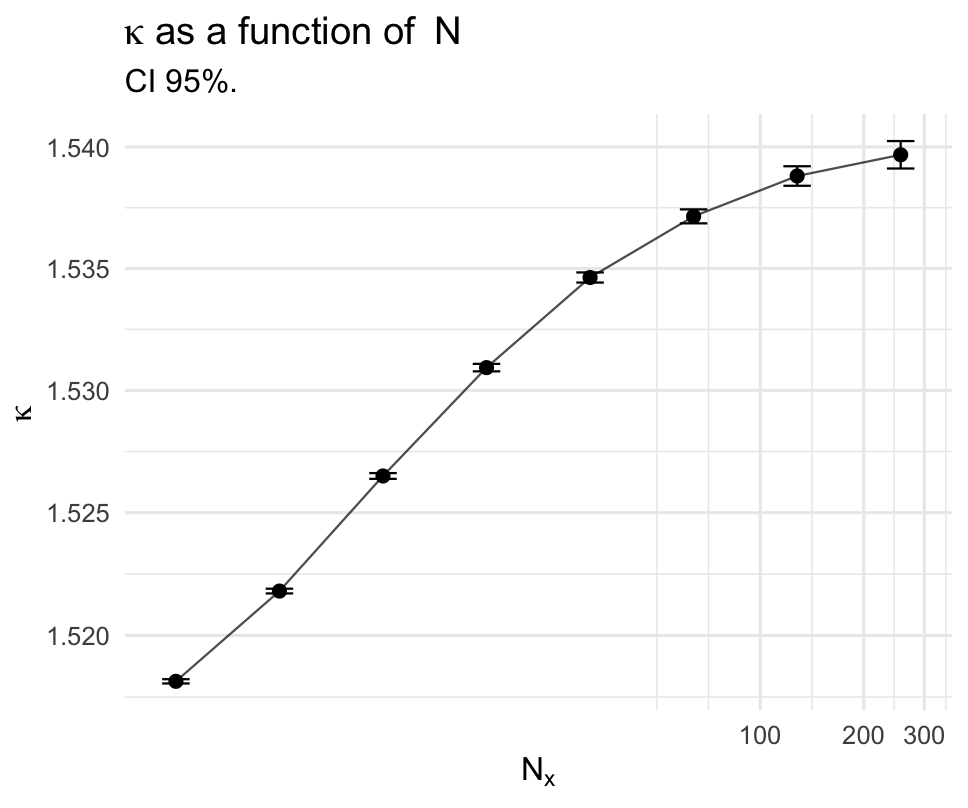}
\caption{Measured values of $\kappa_N$ as a function of the slab length $N$, suggesting convergence toward the predicted limit $\kappa_\infty$. Here
$\kappa_N=\frac{N C_N}{1-C_N}$,
where $C_N$ is the crossing probability of the slab. For each value of $N$, $C_N$ is estimated numerically by sampling random mirror environments, launching a particle from a fixed boundary point, and recording the frequency with which it exits on the opposite side.}
\label{fig:kappa}

\end{center}
\end{figure}


\section{Discussion and Perspectives}

The mirrors model is a deterministic dynamical system in a random environment, and
as such it does not fall into the classical framework of stochastic homogenization.
Nonetheless, the crossing probability $C_N$ exhibits the universal scaling law
\[
   C_N \sim \frac{\kappa_\infty}{N},
\qquad N\to\infty,
\]
with a finite positive conductivity $\kappa_\infty$. 

\noindent Our multiscale analysis identifies the mechanism behind this behavior: on each scale
$2^n$, all path correlations generated by deterministic memory are summarized in a
single correction factor $\Delta_n$. Once $\Delta_n$ is expressed in terms of the
increments $\eta_n(l+1)-\eta_n(l)$, the central feature is the closure assumption, based on increasing independence in one sector of the entry/exit phase space and constant strong correlations on the complement. We have computed the effect of the closure approximation for the two-point correlation
function on the conductivity. This leads to a closed recursion equation for the
conductivity, in which the parameter $\alpha$ collects the effect of trajectory
recollisions. This parameter is computable numerically with high accuracy.

\noindent The limit $\kappa_\infty$ is close to the conductivity constant of the non--backtracking random 
walk; the small difference quantifies the residual memory not captured by a Markovian 
model. This suggests applying the same approach to other deterministic systems, such as
random Lorentz gases with partial density, deterministic cellular automata with local
scatterers, Sinai billiards in finite channels, and transport on random graphs with
local deterministic routing rules.

\noindent A natural next goal is a rigorous control of the error factors $\eta_n(l)$
under suitable geometric decompositions. Establishing exponential decay in $l$ and 
$2^{-n}$ would lead to a fully rigorous proof of normal conductivity in the
three-dimensional mirrors model.

Another natural direction would be to clarify the connection with the recent paper \cite{LefevereTasaki} by investigating how close the matching between incoming and outgoing states in the mirrors model is to the uniform matching rigorously analysed there. Such a comparison could be useful in the search for a rigorous proof of normal conductivity in the mirrors model.

\section*{Data availability}
The data and code supporting the findings of this study are available 

at https://doi.org/10.5281/zenodo.17578332.
\appendix

\section*{Conflict of interest}
The author declares that there is no conflict of interest.

\section{Computation of the Single-Slab Crossing Probability $C_1$}

Let $\gamma=(x_1,\ldots,x_n)$ be a path in $\bbZ^d\times {\mathcal P}$. The probability 
of $\gamma$ is
\begin{equation}
\bbM[\gamma]=\left\langle \prod_{i=1}^n T(x_i,x_{i+1};\pi)\right \rangle.
\label{pathproba0}
\end{equation}
Writing $\un 1=(1,\ldots,1)$, the crossing probability of $\Lambda_N$ can be expressed as
\begin{equation}
C_N=\sum_{\gamma:(\un 1,\e_1)\to O_N^+}\bbM[\gamma],
\label{pathproba}
\end{equation}
where the sum is over paths $\gamma=(x_1,\ldots,x_n)$ with $x_1=(\un 1,\e_1)$ and 
$x_n\in O_N^+$.

For $N=1$ we decompose according to the first step:
\begin{align}
C_1
&=\sum_{\gamma:(\un 1,\e_1)\to O_1^+}\bbM[\gamma]\nonumber\\
&=\left \langle T((\un 1,\e_1),(\un 1+\e_1,\e_1);\pi)\right \rangle\nonumber\\
&\quad+\frac{4}{5}\sum_{\gamma:(\un 1+\e_2,\e_2)\to O_1^+}
\bbM[\gamma\mid\pi(\un 1;\e_1)=\e_2]\nonumber\\
&=\frac{1}{5}+\frac{4}{5}\sum_{\gamma:(\un 1+\e_2,\e_2)\to O_1^+}
\bbM[\gamma\mid\pi(\un 1;\e_1)=\e_2].
\end{align}
Introducing
\begin{equation}
p_\pm = \sum_{\gamma:(\un 1+\e_2,\e_2)\to O_1^\pm}
\bbM[\gamma\mid\pi(\un 1;\e_1)=\e_2],
\label{plpr}
\end{equation}
we have $p_++p_-=1$ because, under the conditioning, there is no loop contained in 
$\Lambda_1$. There is a one-to-one correspondence between paths in $p_+$ and $p_-$, 
except for those exiting at $(\un 1+\e_1,\e_1)$, which are forbidden for $p_-$. Thus
\[
p_+=p_-+\sum_{\gamma:(\un 1+\e_2,\e_2)\to (\un 1+\e_1,\e_1)}
\bbM[\gamma\mid\pi(\un 1;\e_1)=\e_2],
\]
and hence
\[
p_+=\frac{1}{2}(1+q),
\]
with
\begin{equation}
q:=\sum_{\gamma:(\un 1+\e_2,\e_2)\to (\un 1+\e_1,\e_1)}
\bbM[\gamma\mid\pi(\un 1;\e_1)=\e_2].
\end{equation}
We decompose $q=\sum_{n\geq 4}q_n$ according to the number of visited vertices.

For $q_4$ (paths visiting 4 vertices), there are two such paths, each with conditional 
probability $(\frac{1}{5})^3\frac{1}{3}$, so $q_4=(\frac{1}{5})^3\frac{2}{3}$.  
For 6-step loops (visiting 6 vertices) forming rectangles of size $2\times 1$, there are 
six such loops, each with probability $(\frac{1}{5})^5\frac{1}{3}$, giving
$q_6=6 (\frac{1}{5})^5\frac{1}{3}$. Thus
\[
q=0.00597333+\sum_{n\geq 8}q_n,
\]
so
\[
p_+\simeq 0.502987,
\quad
C_1\simeq 0.602389.
\]
Higher $q_n$ can be bounded using the analogy with closed trails on the square lattice; 
each closed orbit of length $n$ has probability at most 
$(\frac{1}{3})^{n/2}(\frac{1}{5})^{n/2-1}$ and there are at most $3^n$ such orbits, so
\[
q_n\leq 5 \left(\frac{3}{5}\right)^{\frac{n}{2}}.
\]


\section{Closure assumption}

Let $\Lambda_{2^n}$ be a slab of width $2^n$, with incoming boundary phase-space sets
$ I^{\pm}$ and outgoing boundary phase-space sets $O^{\pm}$. In this appendix we describe
more explicitly the structural formulas underlying the closure assumption used in the
main text. To avoid confusion with the recursion variables $(x_i,y_i)$ of Sections~4--5,
we denote here by $a_1,a_2$ generic incoming boundary states and by $b_1,b_2$ generic
outgoing boundary states of a single slab of width $2^n$.

We split the set of quadruples $(a_1,b_1,a_2,b_2)$ into five distinct regions on which the
two-point correlation function
\[
\langle t_n(a_1,b_1)\, t_n(a_2,b_2)\rangle
\]
has a different expression. We define the time-reversal operation by
\[
R(\q,\p)=(\q-\p,-\p).
\]
The various Kronecker factors that appear below are direct consequences of the bijective
and reversible character of the mirrors dynamics: they encode the hard constraints
forbidding incompatible pairs of transfers and the forced coincidences coming from the
deterministic structure of the slab-to-slab transfer map.

\begin{figure}[htbp]
\centering

\begin{minipage}{0.32\textwidth}
\centering
\begin{tikzpicture}[scale=0.8, >=Latex]

  \def\Nx{5}
  \def\Ny{4}

  \draw[gray!40] (0,0.5) -- (0,\Ny+0.5);
  \draw[gray!40] (1,0.5) -- (1,\Ny+0.5);
  \draw[gray!40] (\Nx,0.5) -- (\Nx,\Ny+0.5);
  \draw[gray!40] (\Nx+1,0.5) -- (\Nx+1,\Ny+0.5);
  \draw[thick] (1,0.5) rectangle (\Nx,\Ny+0.5);

  \node at (\Nx/2,\Ny/2+0.5) {$\Lambda_{2^n}$};
  \node at (\Nx/2,0.0) {$A^1_{rr}$};

  \fill (1,1) circle (2pt);
  \draw[->,thick] (1,1) -- (1.7,1);
  \node[below left] at (1,1) {$a_1$};

  \fill (0,1) circle (2pt);
  \draw[->,thick] (0,1) -- (-0.7,1);
  \node[above left] at (0,1) {$b_1$};

  \fill (\Nx,3) circle (2pt);
  \draw[->,thick] (\Nx,3) -- (\Nx-0.7,3);
  \node[below left] at (\Nx,3) {$a_2$};

  \fill (\Nx+1,3) circle (2pt);
  \draw[->,thick] (\Nx+1,3) -- (\Nx+1.7,3);
  \node[above left] at (\Nx+1,3) {$b_2$};

  \draw[blue,thick,->]
    (1.7,1)
      .. controls (2.7,0.8) and (2.0,2.6) ..
    (0,1.2);

  \draw[blue,thick,->]
    (\Nx-0.7,3)
      .. controls (3.2,3.9) and (4.2,1.8) ..
    (\Nx+1,3);

\end{tikzpicture}
\end{minipage}
\hfill
\begin{minipage}{0.32\textwidth}
\centering
\begin{tikzpicture}[scale=0.8, >=Latex]

  \def\Nx{5}
  \def\Ny{4}

  \draw[gray!40] (0,0.5) -- (0,\Ny+0.5);
  \draw[gray!40] (1,0.5) -- (1,\Ny+0.5);
  \draw[gray!40] (\Nx,0.5) -- (\Nx,\Ny+0.5);
  \draw[gray!40] (\Nx+1,0.5) -- (\Nx+1,\Ny+0.5);
  \draw[thick] (1,0.5) rectangle (\Nx,\Ny+0.5);

  \node at (\Nx/2,\Ny/2+0.5) {$\Lambda_{2^n}$};
  \node at (\Nx/2,0.0) {$A^2_{rr}$};

  \fill (\Nx,1) circle (2pt);
  \draw[->,thick] (\Nx,1) -- (\Nx-0.7,1);
  \node[below left] at (\Nx,1) {$a_1$};

  \fill (\Nx+1,1) circle (2pt);
  \draw[->,thick] (\Nx+1,1) -- (\Nx+1.7,1);
  \node[above left] at (\Nx+1,1) {$b_1$};

  \fill (\Nx,3) circle (2pt);
  \draw[->,thick] (\Nx,3) -- (\Nx-0.7,3);
  \node[below left] at (\Nx,3) {$a_2$};

  \fill (\Nx+1,3) circle (2pt);
  \draw[->,thick] (\Nx+1,3) -- (\Nx+1.7,3);
  \node[above left] at (\Nx+1,3) {$b_2$};

  \draw[blue,thick,->]
    (\Nx-0.7,1)
      .. controls (3.2,0.5) and (4.0,2.5) ..
    (\Nx+1,1);

  \draw[blue,thick,->]
    (\Nx-0.7,3)
      .. controls (3.2,3.5) and (4.0,1.5) ..
    (\Nx+1,3);

\end{tikzpicture}
\end{minipage}
\hfill
\begin{minipage}{0.32\textwidth}
\centering
\begin{tikzpicture}[scale=0.8, >=Latex]

  \def\Nx{5}
  \def\Ny{4}

  \draw[gray!40] (0,0.5) -- (0,\Ny+0.5);
  \draw[gray!40] (1,0.5) -- (1,\Ny+0.5);
  \draw[gray!40] (\Nx,0.5) -- (\Nx,\Ny+0.5);
  \draw[gray!40] (\Nx+1,0.5) -- (\Nx+1,\Ny+0.5);
  \draw[thick] (1,0.5) rectangle (\Nx,\Ny+0.5);

  \node at (\Nx/2,\Ny/2+0.5) {$\Lambda_{2^n}$};
  \node at (\Nx/2,0.0) {$A_{cr}$};

  \fill (1,1) circle (2pt);
  \draw[->,thick] (1,1) -- (1.7,1);
  \node[below left] at (1,1) {$a_1$};

  \fill (\Nx+1,3) circle (2pt);
  \draw[->,thick] (\Nx+1,3) -- (\Nx+1.7,3);
  \node[above left] at (\Nx+1,3) {$b_1$};

  \fill (\Nx,2) circle (2pt);
  \draw[->,thick] (\Nx,2) -- (\Nx-0.7,2);
  \node[below left] at (\Nx,2) {$a_2$};

  \fill (\Nx+1,1) circle (2pt);
  \draw[->,thick] (\Nx+1,1) -- (\Nx+1.7,1);
  \node[above left] at (\Nx+1,1) {$b_2$};

  \draw[blue,thick,->]
    (1.7,1)
      .. controls (2.5,0.7) and (3.0,3.3) ..
    (\Nx+1,3);

  \draw[blue,thick,->]
    (\Nx-0.7,2)
      .. controls (3.8,2.8) and (4.2,1.2) ..
    (\Nx+1,1);

\end{tikzpicture}
\end{minipage}

\vspace{0.5cm}

\begin{minipage}{0.32\textwidth}
\centering
\begin{tikzpicture}[scale=0.8, >=Latex]

  \def\Nx{5}
  \def\Ny{4}

  \draw[gray!40] (0,0.5) -- (0,\Ny+0.5);
  \draw[gray!40] (1,0.5) -- (1,\Ny+0.5);
  \draw[gray!40] (\Nx,0.5) -- (\Nx,\Ny+0.5);
  \draw[gray!40] (\Nx+1,0.5) -- (\Nx+1,\Ny+0.5);
  \draw[thick] (1,0.5) rectangle (\Nx,\Ny+0.5);

  \node at (\Nx/2,\Ny/2+0.5) {$\Lambda_{2^n}$};
  \node at (\Nx/2,0.0) {$A_{rc}$};

  \fill (1,1) circle (2pt);
  \draw[->,thick] (1,1) -- (1.7,1);
  \node[below left] at (1,1) {$a_2$};

  \fill (\Nx+1,3) circle (2pt);
  \draw[->,thick] (\Nx+1,3) -- (\Nx+1.7,3);
  \node[above left] at (\Nx+1,3) {$b_2$};

  \fill (\Nx,2) circle (2pt);
  \draw[->,thick] (\Nx,2) -- (\Nx-0.7,2);
  \node[below left] at (\Nx,2) {$a_1$};

  \fill (\Nx+1,1) circle (2pt);
  \draw[->,thick] (\Nx+1,1) -- (\Nx+1.7,1);
  \node[above left] at (\Nx+1,1) {$b_1$};

  \draw[blue,thick,->]
    (1.7,1)
      .. controls (2.5,0.7) and (3.0,3.3) ..
    (\Nx+1,3);

  \draw[blue,thick,->]
    (\Nx-0.7,2)
      .. controls (3.8,2.8) and (4.2,1.2) ..
    (\Nx+1,1);

\end{tikzpicture}
\end{minipage}
\hfill
\begin{minipage}{0.32\textwidth}
\centering
\begin{tikzpicture}[scale=0.8, >=Latex]

  \def\Nx{5}
  \def\Ny{4}

  \draw[gray!40] (0,0.5) -- (0,\Ny+0.5);
  \draw[gray!40] (1,0.5) -- (1,\Ny+0.5);
  \draw[gray!40] (\Nx,0.5) -- (\Nx,\Ny+0.5);
  \draw[gray!40] (\Nx+1,0.5) -- (\Nx+1,\Ny+0.5);
  \draw[thick] (1,0.5) rectangle (\Nx,\Ny+0.5);

  \node at (\Nx/2,\Ny/2+0.5) {$\Lambda_{2^n}$};
  \node at (\Nx/2,0.0) {$A_{cc}$};

  \fill (1,1) circle (2pt);
  \draw[->,thick] (1,1) -- (1.7,1);
  \node[below left] at (1,1) {$a_1$};

  \fill (\Nx+1,1) circle (2pt);
  \draw[->,thick] (\Nx+1,1) -- (\Nx+1.7,1);
  \node[above left] at (\Nx+1,1) {$b_1$};

  \fill (1,3) circle (2pt);
  \draw[->,thick] (1,3) -- (1.7,3);
  \node[below left] at (1,3) {$a_2$};

  \fill (\Nx+1,3) circle (2pt);
  \draw[->,thick] (\Nx+1,3) -- (\Nx+1.7,3);
  \node[above left] at (\Nx+1,3) {$b_2$};

  \draw[blue,thick,->]
    (1.7,1)
      .. controls (2.6,0.7) and (3.2,1.5) ..
    (\Nx+1,1);

  \draw[blue,thick,->]
    (1.7,3)
      .. controls (2.6,3.3) and (3.2,2.5) ..
    (\Nx+1,3);

\end{tikzpicture}
\end{minipage}

\caption{Schematic boundary configurations and trajectories corresponding to the
sets $A^1_{rr}$, $A^2_{rr}$, $A_{cr}$, $A_{rc}$, and $A_{cc}$. Each panel shows one
representative configuration belonging to the corresponding class. When a class is
defined as a union of several boundary patterns, the figure displays only one of
these patterns. The points $a_i$ lie on incoming boundaries $I^\pm$ and the points
$b_i$ on outgoing boundaries $O^\pm$. The two trajectories inside $\Lambda_{2^n}$
realize the associated pattern (two returns, one crossing and one return, or two
crossings).}
\label{fig:A-collections}
\end{figure}

The five classes are defined as follows:
\begin{equation}
\begin{aligned}
A^1_{rr}=
\{ a_1\in I^+,\ b_1\in O^-,\ a_2\in I^-,\ b_2\in O^+\}
\cup\{ a_1\in I^-,\ b_1\in O^+,\ a_2\in I^+,\ b_2\in O^-\} .
\end{aligned}
\label{1RR}
\end{equation}

For $A^2_{rr}$, the figure above represents the second branch of the union, namely
the case where both transfers enter from $I^-$ and exit through $O^+$.
\begin{equation}
\begin{aligned}
A^2_{rr}=
\{ a_1\in I^+,\ b_1\in O^-,\ a_2\in I^+,\ b_2\in O^-\}
\cup\{ a_1\in I^-,\ b_1\in O^+,\ a_2\in I^-,\ b_2\in O^+\} .
\end{aligned}
\label{2RR}
\end{equation}

\begin{equation}
\begin{aligned}
A_{cr}=
&\{ a_1\in I^+,\ b_1\in O^+,\ a_2\in I^-,\ b_2\in O^+\}
\cup \{ a_1\in I^+,\ b_1\in O^+,\ a_2\in I^+,\ b_2\in O^-\}
\\
&\cup \{ a_1\in I^-,\ b_1\in O^-,\ a_2\in I^-,\ b_2\in O^+\}
\cup \{ a_1\in I^-,\ b_1\in O^-,\ a_2\in I^+,\ b_2\in O^-\} .
\end{aligned}
\label{CR}
\end{equation}

\begin{equation}
\begin{aligned}
A_{rc}=
&\{ a_2\in I^+,\ b_2\in O^+,\ a_1\in I^-,\ b_1\in O^+\}
\cup \{ a_2\in I^+,\ b_2\in O^+,\ a_1\in I^+,\ b_1\in O^-\}
\\
&\cup \{ a_2\in I^-,\ b_2\in O^-,\ a_1\in I^-,\ b_1\in O^+\}
\cup \{ a_2\in I^-,\ b_2\in O^-,\ a_1\in I^+,\ b_1\in O^-\} .
\end{aligned}
\label{RC}
\end{equation}

\begin{equation}
\begin{aligned}
A_{cc}=
&\{ a_1\in I^+,\ b_1\in O^+,\ a_2\in I^+,\ b_2\in O^+\}
\cup \{ a_1\in I^+,\ b_1\in O^+,\ a_2\in I^-,\ b_2\in O^-\}
\\
&\cup \{ a_1\in I^-,\ b_1\in O^-,\ a_2\in I^+,\ b_2\in O^+\}
\cup \{ a_1\in I^-,\ b_1\in O^-,\ a_2\in I^-,\ b_2\in O^-\} .
\end{aligned}
\label{CC}
\end{equation}

If $(a_1,b_1,a_2,b_2)\in A^1_{rr}$, then
\begin{equation}
\langle t_n(a_1,b_1) t_n(a_2,b_2)\rangle
=c_n(a_1,b_1) c_n(a_2,b_2)\bigl(1+h^1_{n,rr}(a_1,b_1,a_2,b_2)\bigr).
\end{equation}

If $(a_1,b_1,a_2,b_2)\in A^2_{rr}$, then
\begin{align}
\langle t_n(a_1,b_1) t_n(a_2,b_2)\rangle
&=(1-\delta_{a_1a_2})(1-\delta_{b_1b_2})(1-\delta_{a_1Rb_2})(1-\delta_{a_2Rb_1})
\notag\\
&\qquad\cdot c_n(a_1,b_1) c_n(a_2,b_2)
\bigl(1+h^2_{n,rr}(a_1,b_1,a_2,b_2)\bigr)
\notag\\
&\qquad+\bigl(\delta_{a_1a_2}\delta_{b_1b_2}+\delta_{a_1Rb_2}\delta_{a_2Rb_1}\bigr)c_n(a_1,b_1).
\end{align}

If $(a_1,b_1,a_2,b_2)\in A_{cr}$, then
\begin{align}
\langle t_n(a_1,b_1) t_n(a_2,b_2)\rangle
&=(1-\delta_{Rb_1a_2})(1-\delta_{b_1b_2})(1-\delta_{a_1a_2})(1-\delta_{a_1Rb_2})
\notag\\
&\qquad\cdot c_n(a_1,b_1)c_n(a_2,b_2)
\bigl(1+h_{n,cr}(a_1,b_1,a_2,b_2)\bigr).
\end{align}

If $(a_1,b_1,a_2,b_2)\in A_{rc}$, then
\begin{align}
\langle t_n(a_1,b_1) t_n(a_2,b_2)\rangle
&=(1-\delta_{Rb_1a_2})(1-\delta_{b_1b_2})(1-\delta_{a_1a_2})(1-\delta_{a_1Rb_2})
\notag\\
&\qquad\cdot c_n(a_1,b_1)c_n(a_2,b_2)
\bigl(1+h_{n,rc}(a_1,b_1,a_2,b_2)\bigr).
\end{align}

If $(a_1,b_1,a_2,b_2)\in A_{cc}$, then
\begin{align}
\langle t_n(a_1,b_1) t_n(a_2,b_2)\rangle
&=(1-\delta_{b_1b_2})(1-\delta_{a_1a_2})(1-\delta_{Rb_2 a_1})(1-\delta_{Rb_1a_2})
\notag\\
&\qquad\cdot c_n(a_1,b_1)c_n(a_2,b_2)
\bigl(1+h_{n,cc}(a_1,b_1,a_2,b_2)\bigr)
\notag\\
&\qquad+ (\delta_{a_1a_2}\delta_{b_1b_2}+\delta_{a_2Rb_1}\delta_{a_1 Rb_2})c_n(a_1,b_1).
\label{hcc}
\end{align}

All the functions $h_{n,\cdot}$ appearing above are assumed to satisfy the uniform bound
\[
\sup_{a_1,b_1,a_2,b_2}|h_{n,\cdot}(a_1,b_1,a_2,b_2)|\le h(n),
\]
where $h$ is a common positive function such that $h(n)\to 0$ when $n\to\infty$.

The closure assumption used in the main text follows from these structural formulas.
For the $l=2$ terms in the recursion, the relevant two-point functions are initially of
type $A_{cr}$ and $A_{rc}$; in the analysis of $R_2$, the deterministic summation over
outgoing states rewrites the contribution in a form involving the sector $A_{cc}$.
The coefficient $\alpha$ is determined by the limiting value of the rescaled diagonal
contribution $R_{11}$.


\section{Bounds on $R_{11}$, $R_{12}$, and $R_2$}

We summarize uniform estimates on the terms $R_{11}$, $R_{12}$, and $R_2$ defined in 
the main text.
In this appendix, when discussing the terms $R_{11},R_{12},R_2$, we keep the half-slab notation $O_1^\pm$, $O_2^\pm$ of the main text.

\medskip\noindent
\textbf{Lower bound on $R_{11}$.}
$R_{11}$ can be bounded from below uniformly in $n$ by considering trajectories that, 
starting at $y_1$, remain on the left-most slice of $\Lambda_2$ before exiting at 
$x_2\in O_2^-$.

Fix an interface state $y_1\in O_1^+$, identified with the corresponding incoming state for
$\Lambda_2$ by the natural interface correspondence, and let $q_{\mathrm L}$ be its
associated leftmost site.
We bound $(1-c_n)^2 R_{11}$ by summing the squares of the probabilities of disjoint 
very short return events.

Two-step returns (probability $(2d-1)^{-2}$ each): for $j\in\{2,\dots,d\}$ and 
$\sigma\in\{\pm1\}$,
\[
(q_{\mathrm L},\e_1) \to  (q_{\mathrm L}+\sigma \e_j,\sigma \e_j)\to  
(q_{\mathrm L}+\sigma \e_j-\e_1,-\e_1),
\]
exiting at $x_2^{(j,\sigma)}\in O_2^-$. There are $2(d-1)$ such events, and each gives
\[
c_n\bigl(y_1,x_2^{(j,\sigma)}\bigr)\ge (2d-1)^{-2},
\]
so
\[
\sum_{j,\sigma}\bigl(c_n(y_1,x_2^{(j,\sigma)})\bigr)^2\ge
\frac{2(d-1)}{(2d-1)^4}.
\]

Three-step returns (probability $(2d-1)^{-3}$ each): for distinct $j,k\in\{2,\dots,d\}$ 
and $\sigma,\sigma'\in\{\pm1\}$,
\[
(q_{\mathrm L},\e_1) 
\to  (q_{\mathrm L}+\sigma \e_j,\sigma \e_j)
\to  (q_{\mathrm L}+\sigma \e_j+\sigma'\e_k,\sigma'\e_k)
\to (q_{\mathrm L}+\sigma \e_j+\sigma'\e_k-\e_1,-\e_1),
\]
exiting at $x_2^{(j,k,\sigma,\sigma')}\in O_2^-$. There are 
$4(d-1)(d-2)$ such events, giving
\[
\sum_{j\neq k,\sigma,\sigma'}
\bigl(c_n(y_1,x_2^{(j,k,\sigma,\sigma')})\bigr)^2
\ge \frac{4(d-1)(d-2)}{(2d-1)^6}.
\]
Combining we obtain
\[
R_{11}\ge 
\frac{2(d-1)}{(2d-1)^4} + \frac{4(d-1)(d-2)}{(2d-1)^6}.
\]
In $d=3$,
\[
R_{11} \ge \frac{4}{5^4} + \frac{8}{5^6}
\approx 6.912\times 10^{-3}.
\]
\medskip\noindent
\textbf{Bound on $R_{12}$.}

Recall that $R_{12}$ is the contribution to $R_1$ coming from the off-diagonal sector
$y_2\neq y_1$. Writing out all sums explicitly, we have
\begin{equation}
\begin{aligned}
c_n^2(1-c_n)^2R_{12}
=
\sum_{\substack{y_1\in O_1^+,\ y_2\in O_1^+\\ y_2\neq y_1}}
\sum_{x_2\in O_2^-}\sum_{x_3\in O_2^+}
\langle\delta t_n(x_1,y_1)\delta t_n(x_2,y_2)\rangle 
\langle\delta t_n(y_1,x_2)\delta t_n(y_2,x_3)\rangle ,
\end{aligned}
\label{R12explicit}
\end{equation}
with $x_1=(\un 1,\e_1)$ fixed.

For the first two-point function, the pair of transfers belongs to the mixed sector
$A_{cr}$, and since $y_2\neq y_1$ the factor $(1-\delta_{y_1y_2})$ is equal to $1$.
Moreover, as in the main text, the factors $(1-\delta_{x_1x_2})$ and
$(1-\delta_{x_1Ry_2})$ are identically equal to $1$. Hence the closure formula gives
\begin{equation}
\begin{aligned}
\langle t_n(x_1,y_1)t_n(x_2,y_2)\rangle
&=(1-\delta_{Ry_1x_2})\,c_n(x_1,y_1)c_n(x_2,y_2)\bigl(1+O(h(n))\bigr).
\end{aligned}
\label{R12firstclosure}
\end{equation}
Subtracting the product of one-point functions, we obtain
\begin{equation}
\begin{aligned}
\langle\delta t_n(x_1,y_1)\delta t_n(x_2,y_2)\rangle
&=
\Bigl((1-\delta_{Ry_1x_2})(1+O(h(n)))-1\Bigr)c_n(x_1,y_1)c_n(x_2,y_2).
\end{aligned}
\label{R12firstdelta}
\end{equation}
Now, if $\delta_{Ry_1x_2}=1$, then $c_n(y_1,x_2)=0$ by reversibility and bijectivity, so
such terms do not contribute when multiplied by the second two-point function below.
Therefore, on the support relevant to \eqref{R12explicit}, the only surviving part in
\eqref{R12firstdelta} is the closure error, and we may write
\begin{equation}
\langle\delta t_n(x_1,y_1)\delta t_n(x_2,y_2)\rangle
=
O(h(n))\,c_n(x_1,y_1)c_n(x_2,y_2).
\label{R12firstbound}
\end{equation}

We now turn to the second two-point function. It belongs to the sector $A_{rc}$.
Since $y_2\neq y_1$, the factor $(1-\delta_{y_1y_2})$ is equal to $1$, and since
$x_3\in O_2^+$, the factor $(1-\delta_{x_2x_3})$ is also identically equal to $1$.
Moreover, $y_1\in O_1^+$ is an interface state, whereas $Rx_3$ lies on the right incoming
boundary of $\Lambda_2$, so $\delta_{y_1Rx_3}=0$ identically. Hence the closure formula
reduces to
\begin{equation}
\begin{aligned}
\langle t_n(y_1,x_2)t_n(y_2,x_3)\rangle
&=(1-\delta_{y_2Rx_2})\,c_n(y_1,x_2)c_n(y_2,x_3)\bigl(1+O(h(n))\bigr).
\end{aligned}
\label{R12secondclosure}
\end{equation}
Subtracting the product of one-point functions gives
\begin{equation}
\begin{aligned}
\langle\delta t_n(y_1,x_2)\delta t_n(y_2,x_3)\rangle
&=
\Bigl((1-\delta_{y_2Rx_2})(1+O(h(n)))-1\Bigr)c_n(y_1,x_2)c_n(y_2,x_3).
\end{aligned}
\label{R12seconddelta}
\end{equation}
If $\delta_{y_2Rx_2}=1$, then $y_2=Rx_2$, so the transfer $y_2\to x_2$ is the
time-reversal of the trivial interface identification. By bijectivity of the slab
transfer map, this is incompatible with the simultaneous transfer $y_1\to x_2$ unless
$y_1=y_2$, which is excluded in the definition of $R_{12}$. Equivalently, on the
off-diagonal sector $y_2\neq y_1$, such terms do not contribute. Therefore only the
closure error remains and
\begin{equation}
\langle\delta t_n(y_1,x_2)\delta t_n(y_2,x_3)\rangle
=
O(h(n))\,c_n(y_1,x_2)c_n(y_2,x_3).
\label{R12secondbound}
\end{equation}

Inserting \eqref{R12firstbound} and \eqref{R12secondbound} into
\eqref{R12explicit}, we get
\begin{equation}
\begin{aligned}
c_n^2(1-c_n)^2|R_{12}|
&\le
O(h(n)^2)
\sum_{\substack{y_1\in O_1^+,\ y_2\in O_1^+\\ y_2\neq y_1}}
\sum_{x_2\in O_2^-}\sum_{x_3\in O_2^+}
c_n(x_1,y_1)c_n(x_2,y_2)c_n(y_1,x_2)c_n(y_2,x_3).
\end{aligned}
\label{R12bound1}
\end{equation}
Dropping the restriction $y_2\neq y_1$ only enlarges the sum, hence
\begin{equation}
\begin{aligned}
c_n^2(1-c_n)^2|R_{12}|
&\le
O(h(n)^2)
\sum_{y_1\in O_1^+}\sum_{y_2\in O_1^+}\sum_{x_2\in O_2^-}\sum_{x_3\in O_2^+}
c_n(x_1,y_1)c_n(x_2,y_2)c_n(y_1,x_2)c_n(y_2,x_3).
\end{aligned}
\label{R12bound2}
\end{equation}

We now perform the sums successively. First,
\[
\sum_{x_3\in O_2^+} c_n(y_2,x_3)=c_n.
\]
Therefore
\begin{equation}
\begin{aligned}
c_n^2(1-c_n)^2|R_{12}|
&\le
O(h(n)^2)c_n
\sum_{y_1\in O_1^+}\sum_{y_2\in O_1^+}\sum_{x_2\in O_2^-}
c_n(x_1,y_1)c_n(x_2,y_2)c_n(y_1,x_2).
\end{aligned}
\label{R12bound3}
\end{equation}
Next, for fixed $x_2\in O_2^-$,
\[
\sum_{y_2\in O_1^+} c_n(x_2,y_2)=1-c_n,
\]
by the interface identification explained above, since \(x_2\) is viewed as an incoming
state on the right boundary of \(\Lambda_1\), and the exits in \(O_1^+\) are precisely the
complementary events to the crossings into \(O_1^-\). Thus
\begin{equation}
\begin{aligned}
c_n^2(1-c_n)^2|R_{12}|
&\le
O(h(n)^2)c_n(1-c_n)
\sum_{y_1\in O_1^+}\sum_{x_2\in O_2^-}
c_n(x_1,y_1)c_n(y_1,x_2).
\end{aligned}
\label{R12bound4}
\end{equation}
Now, for fixed $y_1\in O_1^+$,
\[
\sum_{x_2\in O_2^-} c_n(y_1,x_2)=1-c_n,
\]
because the transfer through $\Lambda_2$ starting from $y_1$ either returns to $O_2^-$
with probability $1-c_n$ or crosses to $O_2^+$ with probability $c_n$. Hence
\begin{equation}
\begin{aligned}
c_n^2(1-c_n)^2|R_{12}|
&\le
O(h(n)^2)c_n(1-c_n)^2
\sum_{y_1\in O_1^+} c_n(x_1,y_1).
\end{aligned}
\label{R12bound5}
\end{equation}
Finally,
\[
\sum_{y_1\in O_1^+} c_n(x_1,y_1)=c_n.
\]
We conclude that
\[
c_n^2(1-c_n)^2|R_{12}|
\le O(h(n)^2)c_n^2(1-c_n)^2,
\]
and therefore
\begin{equation}
|R_{12}|\le O(h(n)^2).
\label{R12final}
\end{equation}

Thus the off-diagonal part is quadratic in the closure error. Since \(R_{12}=O(h(n)^2)\),
its contribution to \(\Delta_n\) is \(O(c_n h(n)^2)=o(c_n)\) under the sole assumption
\(h(n)\to0\); under the inductive conductive scaling \(c_n\asymp 2^{-n}\), it is therefore
negligible at the conductivity scale.
\medskip\noindent

\textbf{Bound on $R_2$.}

We rewrite the term $R_2$ from the main text with all sums explicit. Recall that
\begin{equation}
\begin{aligned}
c_n^2(1-c_n)^2R_2
&=2\sum_{y_1\in O_1^+}\sum_{x_2\in O_2^-}\sum_{y_2\in O_1^+}\sum_{x_3\in O_2^+}
\langle\delta t_n(x_1,y_1)\delta t_n(x_2,y_2)\rangle \\
&\qquad\qquad\qquad\qquad\qquad\qquad\cdot c_n(y_1,x_2)c_n(y_2,x_3),
\end{aligned}
\label{R2explicit}
\end{equation}
with $x_1=(\un 1,\e_1)$ fixed.

By translation invariance,
\[
\sum_{x_3\in O_2^+} c_n(y_2,x_3)=c_n,
\]
for every $y_2\in O_1^+$. Hence
\begin{equation}
c_n(1-c_n)^2R_2
=
2\sum_{y_1\in O_1^+}\sum_{x_2\in O_2^-}\sum_{y_2\in O_1^+}
\langle\delta t_n(x_1,y_1)\delta t_n(x_2,y_2)\rangle\, c_n(y_1,x_2).
\label{R2afterx3}
\end{equation}

Next, for fixed $x_2\in O_2^-$, the deterministic transfer identity gives
\[
\sum_{y_2\in O_1^+} t_n(x_2,y_2)+\sum_{y_2\in O_1^-} t_n(x_2,y_2)=1,
\]
and therefore also
\[
\sum_{y_2\in O_1^+} \delta t_n(x_2,y_2)+\sum_{y_2\in O_1^-} \delta t_n(x_2,y_2)=0.
\]
Thus
\[
\sum_{y_2\in O_1^+} \delta t_n(x_2,y_2)
=
-\sum_{y_2\in O_1^-} \delta t_n(x_2,y_2),
\]
and \eqref{R2afterx3} becomes
\begin{equation}
c_n(1-c_n)^2R_2
=
-2\sum_{y_1\in O_1^+}\sum_{x_2\in O_2^-}\sum_{y_2\in O_1^-}
\langle\delta t_n(x_1,y_1)\delta t_n(x_2,y_2)\rangle\, c_n(y_1,x_2).
\label{R2rearranged}
\end{equation}

Expanding $\delta t_n=t_n-c_n$, we get
\begin{equation}
\begin{aligned}
c_n(1-c_n)^2R_2
&=
-2\sum_{y_1\in O_1^+}\sum_{x_2\in O_2^-}\sum_{y_2\in O_1^-}
\Bigl(
\langle t_n(x_1,y_1)t_n(x_2,y_2)\rangle
-c_n(x_1,y_1)c_n(x_2,y_2)
\Bigr) \\
&\qquad\qquad\qquad\qquad\qquad\qquad\cdot c_n(y_1,x_2).
\end{aligned}
\label{R2det}
\end{equation}

Before the rearrangement, the relevant two-point function is of mixed type $A_{cr}$.
After replacing the sum over $y_2\in O_1^+$ by the complementary sum over
$y_2\in O_1^-$, the transfer $x_2\to y_2$ becomes, like $x_1\to y_1$, a crossing
transfer through $\Lambda_1$. Hence the two-point function in \eqref{R2det} is now in
the crossing--crossing sector $A_{cc}$, and we may use the structural formula
\eqref{hcc}.
In the present situation, the term $\delta_{a_1a_2}\delta_{b_1b_2}$ in \eqref{hcc}
is identically zero for geometric reasons. The only potentially nonzero forced
contribution is the time-reversed coincidence term
$\delta_{a_2Rb_1}\delta_{a_1Rb_2}c_n(a_1,b_1)$.
In our variables, this term can contribute only if in particular $x_2=Ry_1$.
But then the accompanying factor $c_n(y_1,x_2)$ vanishes by reversibility and
bijectivity. Hence the forced part gives no contribution. Therefore only the closure error term remains, and we obtain
\begin{equation}
\begin{aligned}
c_n(1-c_n)^2R_2
&=
-2\sum_{y_1\in O_1^+}\sum_{x_2\in O_2^-}\sum_{y_2\in O_1^-}
h_{n,cc}(x_1,y_1,x_2,y_2)\,
(1-\delta_{Ry_2x_1})(1-\delta_{Ry_1x_2}) \\
&\qquad\qquad\qquad\qquad\qquad\qquad\cdot
c_n(x_1,y_1)c_n(x_2,y_2)c_n(y_1,x_2).
\end{aligned}
\label{R2hterm}
\end{equation}

Using the uniform bound
\[
|h_{n,cc}(x_1,y_1,x_2,y_2)|\le h(n),
\]
we deduce
\begin{equation}
\begin{aligned}
c_n(1-c_n)^2|R_2|
&\le
2h(n)\sum_{y_1\in O_1^+}\sum_{x_2\in O_2^-}\sum_{y_2\in O_1^-}
c_n(x_1,y_1)c_n(x_2,y_2)c_n(y_1,x_2).
\end{aligned}
\label{R2bound1}
\end{equation}

Now, for fixed $x_2\in O_2^-$, the sum over $y_2\in O_1^-$ is a crossing probability
through $\Lambda_1$, since $x_2$ is on the right side and $y_2$ on the left side. Under the interface identification, a state \(x_2\in O_2^-\) is viewed as an incoming state on the right boundary of \(\Lambda_1\); consequently, exits in \(O_1^-\) correspond to crossings of \(\Lambda_1\), whereas exits in \(O_1^+\) correspond to returns on the same side.
Hence
\[
\sum_{y_2\in O_1^-} c_n(x_2,y_2)=c_n.
\]
Therefore
\begin{equation}
\begin{aligned}
c_n(1-c_n)^2|R_2|
&\le
2h(n)c_n
\sum_{y_1\in O_1^+}\sum_{x_2\in O_2^-}
c_n(x_1,y_1)c_n(y_1,x_2).
\end{aligned}
\label{R2bound2}
\end{equation}

Next, for fixed $y_1\in O_1^+$,
\[
\sum_{x_2\in O_2^-} c_n(y_1,x_2)=1-c_n,
\]
since in $\Lambda_2$ the transfer starting from $y_1$ either crosses to $O_2^+$ with
probability $c_n$ or returns to $O_2^-$ with probability $1-c_n$. Hence
\begin{equation}
\begin{aligned}
c_n(1-c_n)^2|R_2|
&\le
2h(n)c_n(1-c_n)\sum_{y_1\in O_1^+} c_n(x_1,y_1).
\end{aligned}
\label{R2bound3}
\end{equation}

Finally,
\[
\sum_{y_1\in O_1^+} c_n(x_1,y_1)=c_n.
\]
Thus
\[
c_n(1-c_n)^2|R_2|
\le 2h(n)c_n^2(1-c_n),
\]
and therefore
\[
|R_2|\le \frac{2h(n)c_n}{1-c_n}.
\]

Since $c_n\to 0$ as $n\to\infty$, this yields
\[
|R_2|\le c_n O(h(n)).
\]

In particular, the rearrangement of the $y_2$-sum gains one factor $c_n$, because after
the rearrangement the transfer $x_2\to y_2$ is a crossing through $\Lambda_1$. The
remaining sum over $x_2$ produces the factor $1-c_n$ through the transfer $y_1\to x_2$
across $\Lambda_2$.

\end{document}